\documentclass[11pt,reqno]{amsart}
\pdfoutput=1

\usepackage[utopia]{mathdesign}
\usepackage[protrusion=true,expansion=true]{microtype}
\usepackage{setspace}
\usepackage[longnamesfirst]{natbib}
\usepackage{graphicx}
\usepackage[justification=RaggedRight,margin=1em]{subfig}
\usepackage[noend]{algorithmic}

\usepackage[pdftitle={Geometric Computational
  Electrodynamics with Variational Integrators
  and Discrete Differential Forms}, pdfauthor={Ari Stern, Yiying Tong,
  Mathieu Desbrun, and Jerrold E. Marsden}, pdfsubject={math.NA}]{hyperref}

\newtheorem{theorem}{Theorem}[section]

\theoremstyle{definition}

\theoremstyle{remark}

\numberwithin{equation}{section}
\begin{document}

\title[Geometric Computational Electrodynamics]{Geometric
  Computational Electrodynamics\\
  with Variational Integrators\\
  and Discrete Differential Forms}

\author[A.~Stern]{Ari Stern}
\address{Department of Applied and Computational Mathematics\\
California Institute of Technology\\
Pasadena CA 91125}
\curraddr{Department of Mathematics\\
University of California, San Diego\\
La Jolla CA 92093-0112}
\email{astern@math.ucsd.edu}
\thanks{First author's research partially supported by a Gordon and
  Betty Moore Foundation fellowship at Caltech, and by NSF grant
  CCF-0528101.}

\author[Y.~Tong]{Yiying Tong}
\address{Department of Computer Science\\
California Institute of Technology\\
Pasadena CA 91125}
\curraddr{Department of Computer Science and Engineering\\
Michigan State University\\
East Lansing MI 48824}
\email{ytong@msu.edu}

\author[M.~Desbrun]{Mathieu Desbrun}
\address{Department of Computer Science\\
  California Institute of Technology\\
  Pasadena CA 91125} \curraddr{} \email{mathieu@caltech.edu}
\thanks{Second and third authors' research partially supported by NSF
  grants CCR-0133983 and DMS-0453145 and DOE contract
  DE-FG02-04ER25657.}

\author[J.~E.~Marsden]{Jerrold E. Marsden}
\address{Control and Dynamical Systems\\
California Institute of Technology\\
Pasadena CA 91125}
\curraddr{}
\email{jmarsden@caltech.edu}
\thanks{Fourth author's research partially supported by NSF grant CCF-0528101.}


\date{\today}

\dedicatory{}

\begin{abstract}
In this paper, we develop a structure-preserving discretization of the
Lagrangian framework for electromagnetism, combining techniques from
{\em variational integrators} and {\em discrete differential forms}.
This leads to a general family of variational, multisymplectic
numerical methods for solving Maxwell's equations that automatically
preserve key symmetries and invariants.

In doing so, we demonstrate several new results, which apply both to
some well-established numerical methods and to new methods introduced
here.  First, we show that Yee's finite-difference time-domain (FDTD)
scheme, along with a number of related methods, are multisymplectic
and derive from a discrete Lagrangian variational principle. Second,
we generalize the Yee scheme to unstructured meshes, not just in space
but in $4$-dimensional spacetime.  This relaxes the need to take
uniform time steps, or even to have a preferred time coordinate at
all.  Finally, as an example of the type of methods that can be
developed within this general framework, we introduce a new {\em
  asynchronous variational integrator} (AVI) for solving Maxwell's
equations. These results are illustrated with some prototype
simulations that show excellent energy and conservation behavior and
lack of spurious modes, even for an irregular mesh with asynchronous
time stepping.
\end{abstract}

\maketitle



\section{Introduction}

The Yee scheme (also known as finite-difference time-domain, or FDTD)
was introduced in~\citet{Yee1966} and remains one of the most
successful numerical methods used in the field of computational
electromagnetics, particularly in the area of microwave problems.
Although it is not a ``high-order'' method, it is still preferred for
many applications because it preserves important structural features
of Maxwell's equations that other methods fail to capture.  Among
these distinguishing attributes are that the Gauss constraint $ \nabla
\cdot \mathbf{D} = \rho $ is exactly conserved in a discrete sense,
and electrostatic solutions of the form $ \mathbf{E} = - \nabla \phi $
indeed remain stationary in time \citep[see][]{BoRyIn2005}.  In this
paper, we show that these desirable properties are direct consequences
of the variational and discrete differential structure of the Yee
scheme, which mirrors the geometry of Maxwell's equations.  Moreover,
we will show how to construct other variational methods that, as a
result, share these same numerical properties, while at the same time
applying to more general domains.

\subsection{Variational Integrators and Symmetry}
Geometric numerical integrators have been used primarily for the
simulation of classical mechanical systems, where features such as
symplecticity, conservation of momentum, and conservation of energy
are essential.  (For a survey of various methods and applications,
see~\citealp{HaLuWa2006}.)  Among these, {\em variational integrators}
are developed by discretizing the Lagrangian variational principle of
a system, and then requiring that numerical trajectories satisfy a
discrete version of Hamilton's stationary-action principle.  These
methods are automatically symplectic, and they exactly preserve
discrete momenta associated to symmetries of the Lagrangian: for
instance, systems with translational invariance will conserve a
discrete linear momentum, those with rotational invariance will
conserve a discrete angular momentum, etc.  In addition, variational
integrators can be seen to display good long-time energy behavior,
without artificial numerical damping \citep[see][for a comprehensive
overview of key results]{MaWe2001}.

This variational approach was extended to discretizing general
multisymplectic field theories, with an application to nonlinear wave
equations, in~\citet{MaPaSh1998} and~\citet{MaPeShWe2001}, which
developed the multisymplectic approach for continuum
mechanics. Building on this work, \citet{LeMaOrWe2003} introduced {\em
  asynchronous variational integrators} (AVIs), with which it becomes
possible to choose a different time step size for each element of the
spatial mesh, while still preserving the same variational and
geometric structure as uniform-time-stepping schemes.  These methods
were implemented and shown to be not only practical, but in many cases
superior to existing methods for problems such as nonlinear
elastodynamics. Some further developments are given
in~\citet{LeMaOrWe2004}.

While there have been attempts to apply the existing AVI theory to
computational electromagnetics, these efforts encountered a
fundamental obstacle.  The key symmetry of Maxwell's equations is not
rotational or translational symmetry, as in mechanics, but a {\em
  differential gauge symmetry}. Without taking additional care to
preserve this gauge structure, even variational integrators cannot be
expected to capture the geometry of Maxwell's equations.  As will be
explained, we overcome this obstacle by combining variational methods
with discrete differential forms and operators.  This
differential/gauge structure also turns out to be important for the
numerical performance of the method, and is one of the hallmarks of
the Yee scheme.

\subsection{Preserving Discrete Differential Structure}
As motivation, consider the basic relation $ \mathbf{B} = \nabla
\times \mathbf{A} $, where $\mathbf{B}$ is the magnetic flux and
$\mathbf{A}$ is the magnetic vector potential.  Because of the vector
calculus identities $ \nabla \cdot \nabla \times = 0 $ and $ \nabla
\times \nabla = 0 $, this equation has two immediate and important
consequences.  First, $\mathbf{B}$ is automatically divergence-free.
Second, any transformation $ \mathbf{A} \mapsto \mathbf{A} + \nabla f
$ has no effect on $\mathbf{B}$; this describes a gauge symmetry, for
which the associated conserved momentum is $ \nabla \cdot \mathbf{D} -
\rho $ (which must be zero by Gauss' law).  A similar argument also
explains the invariance of electrostatic solutions, since $ \mathbf{E}
= - \nabla \phi $ is curl-free and invariant under constant shifts in
the scalar potential $\phi$.  Therefore, a proper variational
integrator for electromagnetism should also preserve a discrete analog
of these differential identities.

This can be done by viewing the objects of electromagnetism not as
vector fields, but as {\em differential forms} in $4$-dimensional
spacetime, as is typically done in the literature on classical field
theory.  Using a discrete exterior calculus (called DEC) as the
framework to discretize these differential forms, we find that the
resulting variational integrators automatically respect discrete
differential identities such as $ \mathrm{d} ^2 = 0 $ (which
encapsulates the previous div-curl-grad relations) and Stokes'
theorem.  Consequently, they also respect the gauge symmetry of
Maxwell's equations, and therefore preserve the associated discrete
momentum.

\subsection{Geometry has Numerical Consequences}
The Yee scheme, as we will show, is a method of precisely this type,
which gives a new explanation for many of its previously observed {\em
  a posteriori} numerical qualities.  For instance, one of its notable
features is that the electric field $\mathbf{E}$ and magnetic field
$\mathbf{H}$ do not live at the same discrete space or time locations,
but at separate nodes on a staggered lattice.  The reason why this
particular setup leads to improved numerics is not obvious: if we view
$\mathbf{E}$ and $\mathbf{H}$ simply as vector fields in
$3$-space---the exact same type of mathematical object---why shouldn't
they live at the same points?  Indeed, many finite element method
(FEM) approaches do exactly this, resulting in a ``nodal''
discretization.  However, from the perspective of differential forms
in spacetime, it becomes clear that the staggered-grid approach is
more faithful to the structure of Maxwell's equations: as we will see,
$\mathbf{E}$ and $\mathbf{H}$ come from objects that are dual to one
another (the spacetime forms $F$ and $ G = *F $), and hence they
naturally live on two staggered, dual meshes.

The argument for this approach is not merely a matter of theoretical
interest: the geometry of Maxwell's equations has important practical
implications for numerical performance.  For instance, the
vector-field-based discretization, used in nodal FEM, results in
spurious $3$-D artifacts due to its failure to respect the underlying
geometric structure.  The Yee scheme, on the other hand, produces
resonance spectra in agreement with theory, without spurious modes
\citep[see][]{BoRyIn2005}.  Furthermore, it has been shown
in~\citet{HaAs2001} that staggered-grid methods can be used to develop
fast numerical methods for electromagnetism, even for problems in
heterogeneous media with highly discontinuous material parameters such
as conductivity and permeability.

By developing a structure-preserving, geometric discretization of
Maxwell's equations, not only can we better understand the Yee scheme
and its characteristic advantages, but we can also construct more
general methods that share its desirable properties.  This family of
methods includes the ``Yee-like'' scheme of~\citet{BoKe2000}, which
presented the first extension of Yee's scheme to unstructured grids
(e.g., simplicial meshes rather than rectangular lattices).  General
methods like these are highly desirable: rectangular meshes are not
always practical or appropriate to use in applications where domains
with curved and oblique boundaries are needed \citep[see, for
instance][]{ClWe2002}.  By allowing general discretizations while
still preserving geometry, one can combine the best attributes of the
FEM and Yee schemes.

\subsection{Contributions}
Using DEC as a structure-preserving, geometric framework for general
discrete meshes, we have obtained the following results:
\begin{enumerate}
\item The Yee scheme is actually a variational integrator: that is, it
  can be obtained by applying Hamilton's principle of stationary
  action to a discrete Lagrangian.
\item Consequently, the Yee scheme is \emph{multisymplectic} and
  \emph{preserves discrete momentum maps} (i.e., conserved quantities
  analogous to the continuous case of electromagnetism).  In
  particular, the Gauss constraint is understood as a discrete
  momentum map of this integrator, while the preservation of
  electrostatic potential solutions corresponds to the identity $
  \mathrm{d} ^2 = 0 $, where $\mathrm{d}$ is the discrete exterior
  derivative operator.
\item We also create a foundation for more general schemes, allowing
  {\em arbitrary discretizations of spacetime}, not just uniform time
  steps on a spatial mesh.  One such scheme, introduced here, is a new
  asynchronous variational integrator (AVI) for Maxwell's equations,
  where each spatial element is assigned its own time step size and
  evolves ``asynchronously'' with its neighbors.  This means that one
  can choose to take small steps where greater refinement is needed,
  while still using larger steps for other elements.  Since refining
  one part of the mesh does not restrict the time steps taken
  elsewhere, an AVI can be computationally efficient and numerically
  stable with fewer total iterations.  In addition to the AVI scheme,
  we briefly sketch how completely covariant spacetime integrators for
  electromagnetism can be implemented, without even requiring a 3+1
  split into space and time components.
\end{enumerate}

\subsection{Outline}
We will begin by reviewing Maxwell's equations: first developing the
differential forms expression from a Lagrangian variational principle,
and next showing how this is equivalent to the familiar vector
calculus formulation.  We will then motivate the use of DEC for
computational electromagnetics, explaining how electromagnetic
quantities can be modeled using discrete differential forms and
operators on a spacetime mesh.  These DEC tools will then be used to
set up the discrete Maxwell's equations, and to show that the
resulting numerical algorithm yields the Yee and Bossavit--Kettunen
schemes as special cases, as well as a new AVI method.  Finally, we
will demonstrate that the discrete Maxwell's equations can also be
derived from a discrete variational principle, and will explore its
other discrete geometric properties, including multisymplecticity and
momentum map preservation.

\section{Maxwell's Equations}\label{MaxwellEquations}
This section quickly reviews the differential forms approach to
electromagnetism, in preparation for the associated discrete
formulation given in the next section. For more details, the reader
can refer to~\citet{Bossavit1998} and~\citet{GrKo2004}.

\subsection{From Vector Fields to Differential Forms}
Maxwell's equations, without free sources of charge or current, are
traditionally expressed in terms of four vector fields in $3$-space:
the electric field $\mathbf{E}$, magnetic field $\mathbf{H}$, electric
flux density $\mathbf{D}$, and magnetic flux density $\mathbf{B}$.  To
translate these into the language of differential forms, we begin by
replacing the electric field with a $1$-form $E$ and the magnetic flux
density by a $2$-form $B$.  These have the coordinate expressions
\begin{align*} 
  E &= E _x \;\mathrm{d}x + E _y \;\mathrm{d}y + E _z \;\mathrm{d}z \\
  B &= B _x \;\mathrm{d}y \wedge \mathrm{d} z + B _y \;\mathrm{d}z
  \wedge \mathrm{d} x + B _z \;\mathrm{d}x \wedge \mathrm{d} y,
\end{align*} 
where $\mathbf{E} = (E _x, E _y, E _z) $ and $ \mathbf{B} = (B _x , B
_y, B _z) $.  The motivation for choosing $E$ as a $1$-form and $B$ as
a $2$-form comes from the integral formulation of Faraday's law,
\begin{equation*} 
  \oint _C \mathbf{E} \cdot \mathrm{d} \mathbf{l} =
  -\frac{\mathrm{d}}{\mathrm{d}t} \int _S \mathbf{B} \cdot \mathrm{d}
  \mathbf{A},
\end{equation*} 
where $\mathbf{E}$ is integrated over curves and $\mathbf{B}$ is
integrated over surfaces.  Similarly, Amp\`ere's law,
\begin{equation*} 
  \oint _C \mathbf{H} \cdot \mathrm{d} \mathbf{l} =
  \frac{\mathrm{d}}{\mathrm{d}t} \int _S
  \mathbf{D} \cdot \mathrm{d} \mathbf{A},
\end{equation*} 
integrates $\mathbf{H}$ over curves and $\mathbf{D}$ over surfaces, so
we can likewise introduce a $1$-form $H$ and a $2$-form $D$.

Now, $\mathbf{E}$ and $\mathbf{B}$ are related to $\mathbf{D}$ and
$\mathbf{H}$ through the usual constitutive relations
\begin{equation*} 
  \mathbf{D} = \epsilon \mathbf{E}, \qquad \mathbf{B} = \mu \mathbf{H} .
\end{equation*} 
As shown in~\citet{BoKe2000}, we can view $\epsilon$ and $\mu$ as
corresponding to Hodge operators $*_\epsilon $ and $ * _\mu $, which
map the $1$-form ``fields'' to $2$-form ``fluxes'' in space.
Therefore, this is compatible with viewing $E$ and $H$ as $1$-forms,
and $D$ and $B$ as $2$-forms.

Note that in a vacuum, with $ \epsilon = \epsilon _0 $ and $ \mu = \mu
_0 $ constant, one can simply express the equations in terms of
$\mathbf{E}$ and $\mathbf{B}$, choosing appropriate geometrized units
such that $\epsilon _0 = \mu _0 = c = 1 $, and hence ignoring the
distinction between $\mathbf{E}$ and $\mathbf{D}$ and between
$\mathbf{B}$ and $\mathbf{H}$.  This is typically the most familiar
form of Maxwell's equations, and the one that most students of
electromagnetism first encounter.  In this presentation, we will
restrict ourselves to the vacuum case with geometrized units; for
geometric clarity, however, we will always distinguish between the
$1$-forms $E$ and $H$ and the $2$-forms $D$ and $B$.

Finally, we can incorporate free sources of charge and current by
introducing the {\em charge density $3$-form} $ \rho \;\mathrm{d}x
\wedge \mathrm{d} y \wedge \mathrm{d} z $, as well as the {\em current
  density $2$-form} $ J = J _x \;\mathrm{d}y \wedge \mathrm{d} z + J
_y \;\mathrm{d}z \wedge \mathrm{d} x + J _z \;\mathrm{d}x \wedge
\mathrm{d} y $.  These are required to satisfy the continuity of
charge condition $ \partial _t \rho + \mathrm{d} J = 0 $, which can be
understood as a conservation law (in the finite volume sense).

\subsection{The Faraday and Maxwell $2$-Forms}
In Lorentzian spacetime, we can now combine $E$ and $B$ into a single
object, the Faraday $2$-form
\begin{equation*} 
  F = E \wedge \mathrm{d}  t + B.
\end{equation*} 
There is a theoretical advantage to combining the electric field and
magnetic flux into a single spacetime object: this way,
electromagnetic phenomena can be described in a relativistically
covariant way, without favoring a particular split of spacetime into
space and time components.  In fact, we can turn the previous
construction around: take $F$ to be the fundamental object, with $E$
and $B$ only emerging when we choose a particular coordinate frame.
Taking the Hodge star of $F$, we also get a dual $2$-form
\begin{equation*} 
  G = *F = H \wedge \mathrm{d}  t - D ,
\end{equation*} 
called the Maxwell $2$-form.  The equation $ G = *F $ describes the
dual relationship between $E$ and $B$ on one hand, and $D$ and $H$ on
the other, that is expressed in the constitutive relations.

\subsection{The Source $3$-Form} Likewise, the charge density $\rho$
and current density $J$ can be combined into a single spacetime
object, the source $3$-form
\begin{equation*}
  \mathcal{J} = J \wedge \mathrm{d} t - \rho .
\end{equation*}
Having defined $\mathcal{J}$ in this way, the continuity of charge
condition simply requires that $\mathcal{J}$ be closed, i.e., $
\mathrm{d} \mathcal{J} = 0 $.

\subsection{Electromagnetic Variational Principle}
Let $A$ be the electromagnetic potential $1$-form, satisfying $ F =
\mathrm{d} A $, over the spacetime manifold $X$.  Then define the
$4$-form Lagrangian density
\begin{equation*} 
  \mathcal{L} = - \frac{1}{2} \mathrm{d}  A \wedge *\mathrm{d}  A + A \wedge \mathcal{J} ,
\end{equation*} 
and its associated action functional
\begin{equation*} 
  S[A] = \int _X \mathcal{L}.
\end{equation*} 
Now, take a variation $\alpha$ of $A$, where $\alpha$ vanishes on the
boundary $ \partial X $.  Then the variation of the action functional
along $\alpha$ is
\begin{align*} 
  \mathbf{d} S [A] \cdot \alpha & = \left. \frac{\mathrm{d}
    }{\mathrm{d} \epsilon } \right\rvert _{ \epsilon = 0 } S[A + \epsilon
  \alpha ] \\
  &= \int _X \left( - \mathrm{d} \alpha \wedge *\mathrm{d} A + \alpha
    \wedge \mathcal{J} \right) \\
  &= \int _X \alpha \wedge \left( - \mathrm{d} {*\mathrm{d}A} +
    \mathcal{J} \right) ,
\end{align*} 
where in this last equality we have integrated by parts, using the
fact that $\alpha$ vanishes on the boundary.  Hamilton's principle of
stationary action requires this variation to be equal to zero for
arbitrary $\alpha$, thus implying the electromagnetic Euler--Lagrange
equation,
\begin{equation}
  \mathrm{d} {*\mathrm{d} A} = \mathcal{J} . \label{maxa}
\end{equation}

\subsection{Variational Derivation of Maxwell's Equations}
Since $ G = *F = *\mathrm{d} A $, then clearly \autoref{maxa} is
equivalent to $ \mathrm{d} G = \mathcal{J} $.  Furthermore, since $
\mathrm{d} ^2 = 0 $, it follows that $ \mathrm{d} F = \mathrm{d} ^2 A
= 0 $.  Hence, Maxwell's equations with respect to the Maxwell and
Faraday $2$-forms can be written as
\begin{align}
  \mathrm{d}  F & =   0  \label{maxf1}\\
  \mathrm{d} G & = \mathcal{J} \label{maxf2}
\end{align}

Suppose now we choose the standard coordinate system $ (x, y, z, t) $
on Minkowski space $X = \mathbb{R}^{3,1}$, and define $E$ and $B$
through the relation $ F = E \wedge \mathrm{d} t + B $.  Then a
straightforward calculation shows that \autoref{maxf1} is equivalent
to
\begin{align}
  \nabla \times \mathbf{E} + \partial _t  \mathbf{B} & = 0 \label{curle} \\
  \nabla \cdot \mathbf{B} &= 0 \label{divb}.
\end{align}
Likewise, if $G = *F = H \wedge \mathrm{d} t - D $, then
\autoref{maxf2} is equivalent to
\begin{align}
  \nabla \times \mathbf{H} - \partial _t \mathbf{D} & = \mathbf{J} \label{curlh} \\
  \nabla \cdot \mathbf{D} & = \rho \label{divd}.
\end{align}
Hence this Lagrangian, differential forms approach to Maxwell's
equations is \emph{strictly equivalent} to the more classical vector
calculus formulation in smooth spacetime.  However, in discrete
spacetime, we will see that the differential forms version is
\emph{not} equivalent to an arbitrary vector field discretization,
but rather implies a particular choice of discrete objects.

\subsection{Generalized Hamilton--Pontryagin Principle for Maxwell's
  Equations}
We can also derive Maxwell's equations by using a mixed variational
principle, similar to the Hamilton--Pontryagin principle introduced by
\citet{YoMa2006b} for classical Lagrangian mechanics.  To do this, we
treat $A$ and $F$ as separate fields, while $G$ acts as a Lagrange
multiplier, weakly enforcing the constraint $ F = \mathrm{d} A $.
Define the extended action to be
\begin{equation*}
  S[A,F,G] = \int _X \left[  - \frac{1}{2} F \wedge *F + A \wedge
    \mathcal{J} + \left( F - \mathrm{d} A \right) \wedge G \right] .
\end{equation*}
Then, taking the variation of the action along some $ \alpha, \phi ,
\gamma $ (vanishing on $ \partial X $), we have
\begin{align*}
  \mathbf{d} S [A,F,G] \cdot \left( \alpha, \phi , \gamma \right) &=
  \int _X \left[ - \phi \wedge * F + \alpha \wedge \mathcal{J} +
    \left( \phi - \mathrm{d} \alpha \right) \wedge G + \left( F -
      \mathrm{d} A \right) \wedge \gamma \right] \\
  &= \int _X \left[ \alpha \wedge \left( \mathcal{J} - \mathrm{d} G
    \right) + \phi \wedge \left( G - *F \right) + \left( F -
      \mathrm{d} A \right) \wedge \gamma \right] .
\end{align*}
Therefore, setting this equal to zero, we get the equations
\begin{equation*}
  \mathrm{d} G = \mathcal{J} , \qquad G = *F, \qquad F = \mathrm{d} A .
\end{equation*}
This is precisely equivalent to Maxwell's equations, as derived above.
However, this approach provides some additional insight into the
geometric structure of electromagnetics: the gauge condition $ F =
\mathrm{d} A $ and constitutive relations $ G = *F $ are {\em
  explicitly included in the equations of motion}, as a direct result
of the variational principle.

\subsection{Reducing the Equations}
When solving an initial value problem, it is not necessary to use all
of Maxwell's equations to evolve the system forward in time.  In fact,
the curl equations \eqref{curle} and \eqref{curlh} automatically
conserve the quantities $ \nabla \cdot \mathbf{B} $ and $
\nabla \cdot \mathbf{D} - \rho $.  Therefore, the divergence
equations \eqref{divb} and \eqref{divd} can be viewed simply as
constraints on initial conditions, while the curl equations completely
describe the time evolution of the system.

There are a number of ways to see why we can justify eliminating the
divergence equations.  A straightforward way is to take the divergence
of equations \eqref{curle} and \eqref{curlh}.  Since $ \nabla \cdot
\nabla \times = 0 $, we are left with
\begin{equation*} 
  \partial _t  \left( \nabla \cdot \mathbf{B}\right)  = 0,
  \qquad  \partial _t  \left( \nabla \cdot \mathbf{D} \right)  
  + \nabla \cdot \mathbf{J} = \partial _t \left( \nabla \cdot
    \mathbf{D} - \rho \right) = 0  .
\end{equation*} 
Therefore, if the divergence constraints are satisfied at the initial
time, then they are satisfied for all time, since the divergence terms
are constant.

Another approach is to notice that Maxwell's equations depend only on
the exterior derivative $ \mathrm{d} A $ of the electromagnetic
potential, and not on the value of $A$ itself.  Therefore, the system
has a \emph{gauge symmetry}: any gauge transformation $ A \mapsto A +
\mathrm{d} f $ leaves $ \mathrm{d} A $, and hence Maxwell's equations,
unchanged.  Choosing a time coordinate, we can then partially fix the
gauge so that the electric scalar potential $ \phi = A \left( \partial
  / \partial t \right) = 0$ (the so-called Weyl gauge or temporal
gauge), and so $A$ has only spatial components.  In fact, these three
remaining components correspond to those of the usual vector potential
$\mathbf{A}$.  The reduced Euler--Lagrange equations in this gauge
consist only of \autoref{curlh}, while the remaining gauge symmetry $
\mathbf{A} \mapsto \mathbf{A} + \nabla f $ yields a momentum map that
automatically preserves $ \nabla \cdot \mathbf{D} - \rho $ in time.
Equations \eqref{curle} and \eqref{divb} are automatically preserved
by the identity $ \mathrm{d} ^2 A = 0 $; they are not actually part of
the Euler--Lagrange equations. A more detailed exposition of these
calculations will be given in \autoref{section:gauge}.

\section{Discrete Forms in Computational Electromagnetics}

In this section, we give a quick review of the fundamental objects and
operations of Discrete Exterior Calculus (DEC), a structure-preserving
calculus of discrete differential forms.  By construction, DEC
automatically preserves a number of important geometric structures,
and hence it provides a fully discrete analog of the tools used in the
previous section to express the differential forms version of
Maxwell's equations.  In subsequent sections, we will use this
framework to formulate Maxwell's equations discretely, emulating the
continuous version.

\subsection{Rationale Behind DEC for Computational Electromagnetics}
Modern computational electromagnetism started in the 1960s, when the
finite element method (FEM), based on \emph{nodal basis functions},
was used successfully to discretize the differential equations
governing $2$-{D} static problems formulated in terms of a scalar
potential. Unfortunately, the initial success of the FEM approach
appeared unable to carry over to $3$-{D} problems without spurious
numerical artifacts.  With the introduction of \emph{edge elements}
in~\citet{Nedelec1980} came the realization that a better
discretization of the geometric structure of Maxwell's electromagnetic
theory was key to overcoming this obstacle (see~\citealp{GrKo2004} for
more historical details). Mathematical tools developed by Weyl and
Whitney in the 1950s, in the context of algebraic topology, turned out
to provide the necessary foundations on which robust numerical
techniques for electromagnetism can be built, as detailed
in~\citet{Bossavit1998}.

\subsection{Discrete Differential Forms and Operators}
\label{sec:dec}

In this section, we show how to define differential forms and
operators on a discrete mesh, in preparation to use this framework for
computational modeling of classical fields.  By construction, the
calculus of discrete differential forms automatically preserves a
number of important geometric structures, including Stokes' theorem,
integration by parts (with a proper treatment of boundaries), the
de~Rham complex, Poincar\'{e} duality, Poincar\'{e}'s lemma, and Hodge
theory.  Therefore, this provides a suitable foundation for the
coordinate-free discretization of geometric field theories.  In
subsequent chapters, we will also use these discrete differential
forms as the space of fields on which we will define discrete
Lagrangian variational principles.

The particular ``flavor'' of discrete differential forms and operators
we will be using is known as {\em discrete exterior calculus}, or {\em
  DEC} for short; see \citet{Hirani2003,Leok2004}.  (For related
efforts in this direction, see also~\citealp{Harrison2005}
and~\citealp{ArFaWi2006}.)  Guided by Cartan's exterior calculus of
differential forms on smooth manifolds, DEC is a discrete calculus
developed, \emph{ab initio}, on discrete manifolds, so as to maintain
the covariant nature of the quantities involved.  This computational
tool is based on the notion of discrete chains and cochains, used as
basic building blocks for compatible discretizations of important
geometric structures such as the de~Rham complex
\citep{DeKaTo2008}. The chain and cochain representations are not only
attractive from a computational perspective due to their conceptual
simplicity and elegance; as we will see, they also originate from a
theoretical framework defined by~\citet{Whitney1957}, who introduced
the Whitney and de~Rham maps that establish an isomorphism between
simplicial cochains and Lipschitz differential forms.

\subsubsection{Mesh and Dual Mesh} DEC is concerned with problems in
which the smooth $n$-dimensional manifold $X$ is replaced by a
discrete mesh---precisely, by a cell complex that is manifold, admits
a metric, and is orientable.  The simplest example of such a mesh is a
finite simplicial complex, such as a triangulation of a
$2$-dimensional surface.  We will generally denote the complex by $K$,
and a cell in the complex by $ \sigma $.

\begin{figure}
\centerline{\includegraphics[width=0.7\linewidth]{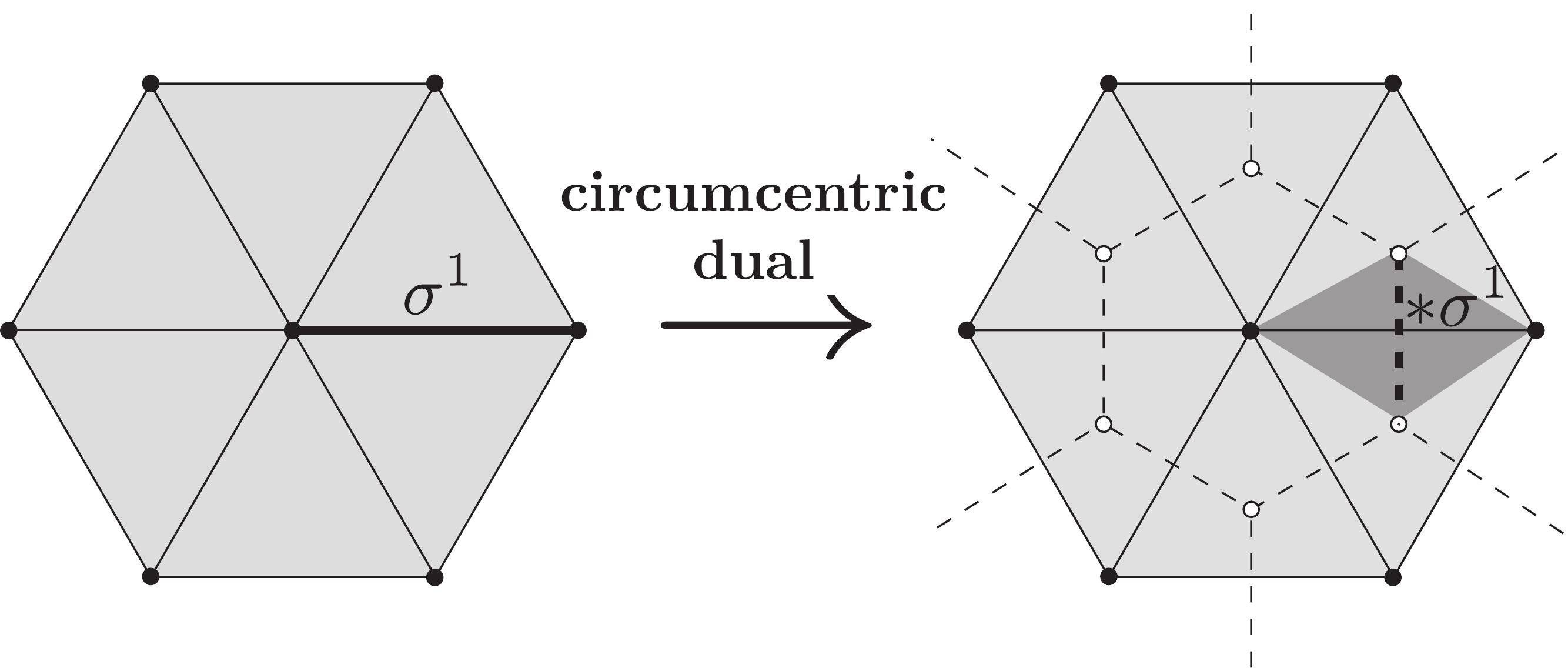}}
\caption[primal-dual]{Given a $2$-D simplicial mesh (left), we can
  construct its circumcentric dual mesh, called the Voronoi diagram of
  the primal mesh (right).  In bold, we show one particular primal
  edge $ \sigma ^1 $ (left) and its corresponding dual edge $ * \sigma
  ^1 $ (right); the convex hull of these cells $
  \operatorname{CH}(\sigma ^1 , *\sigma ^1 ) $ is shaded dark
  grey.}  \label{fig:primal-dual}
\end{figure}

Given a mesh $K$, one can construct a \emph{dual mesh} $ *K $, where
each $k$-cell $\sigma$ corresponds to a dual $(n-k)$-cell $ *\sigma
$. ($ *K $ is ``dual'' to $K $ in the sense of a graph dual.)  One way
to do this is as follows: place a dual vertex at the circumcenter of
each $n$-simplex, then connect two dual vertices by an edge wherever
the corresponding $n$-simplices share an $(n-1)$-simplex, and so on.
This is called the \emph{circumcentric dual}, and it has the important
property that primal and dual cells are automatically orthogonal to
one another, which is advantageous when defining an inner product (as
we will see later in this section). For example, the circumcentric
dual of a Delaunay triangulation, with the Euclidean metric, is its
corresponding Voronoi diagram (see \autoref{fig:primal-dual}).  For
more on the dual relationship between Delaunay triangulations and
Voronoi diagrams, a standard reference is~\citet{ORourke1998}.  A
similar construction of the circumcenter can be carried out for
higher-dimensional Euclidean simplicial complexes, as well as for
simplicial meshes in Minkowski space.  Note that, in both the
Euclidean and Lorentzian cases, the circumcenter may actually lie {\em
  outside} the simplex if it has a very bad aspect ratio, underscoring
the importance of mesh quality for good numerical results.

There are alternative ways to define the dual mesh---for example,
placing dual vertices at the barycenter rather than the
circumcenter---but we will use the circumcentric dual unless otherwise
noted. Note that a refined definition of the dual mesh, where dual
cells at the boundary are restricted to $K$, will be discussed
in~\autoref{sec:boundaryDEC} to allow proper enforcement of boundary
conditions in computational electromagnetics.

\subsubsection{Discrete Differential Forms}
The fundamental objects of DEC are discrete differential forms.  A
discrete $k$-form $\alpha^k$ assigns a real number to each oriented
$k$-dimensional cell $ \sigma ^k $ in the mesh $K$.  (The superscripts
$k$ are not actually required by the notation, but they are often
useful as reminders of what order of form or cell we are dealing
with.)  This value is denoted by $ \left\langle \alpha ^k , \sigma ^k
\right\rangle $, and can be thought of as the value of $ \alpha ^k $
``integrated over'' the element $ \sigma ^k $, i.e.,
\begin{equation*} 
  \left\langle \alpha , \sigma \right\rangle \equiv \int _{ \sigma }
  \alpha .
\end{equation*} 
For example, $0$-forms assign values to vertices, $1$-forms assign
values to edges, etc.  We can extend this to integrate over discrete
paths by linearity: simply add the form's values on each cell in the
path, taking care to flip the sign if the path is oriented opposite
the cell.  Formally, these ``paths'' of $k$-dimensional elements are
called \emph{chains}, and discrete differential forms are
\emph{cochains}, where $ \left\langle \cdot, \cdot \right\rangle $ is
the pairing between cochains and chains.

Differential forms can be defined either on the mesh $K$ or on its
dual $ *K $.  We will refer to these as \emph{primal forms} and
\emph{dual forms} respectively.  Note that there is a natural
correspondence between primal $k$-forms and dual $(n-k)$-forms, since
each primal $k$-cell has a dual $(n-k)$-cell.  This is an important
property that will be used below to define the discrete Hodge star
operator.

\subsubsection{Exterior Derivative}
The discrete exterior derivative $\mathrm{d}$ is constructed to
satisfy Stokes' theorem, which in the continuous sense is written
\begin{equation*} 
  \int _\sigma \mathrm{d} \alpha = \int _{ \partial \sigma } \alpha .
\end{equation*} 
Therefore, if $\alpha$ is a discrete differential $k$-form, then the
$( k + 1 ) $-form $ \mathrm{d} \alpha $ is defined on any $( k + 1 )
$-chain $\sigma$ by
\begin{equation*} 
  \left\langle \mathrm{d}  \alpha , \sigma \right\rangle =
  \left\langle \alpha , \partial \sigma \right\rangle ,
\end{equation*} 
where $\partial \sigma $ is the $ k $-chain boundary of $\sigma $.
For this reason, $\mathrm{d}$ is often called the \emph{coboundary}
operator in cohomology theory.

\subsubsection{Diagonal Hodge Star}
The discrete Hodge star transforms $k$-forms on the primal mesh into
$(n-k)$-forms on the dual mesh, and vice-versa.  In our setup, we will
use the so-called diagonal (or mass-lumped) approximation of the Hodge
star \citep{Bossavit1998} because of its simplicity, but note that
higher-order accurate versions can be substituted.  Given a discrete
form $\alpha$, its \emph{Hodge star} $*\alpha$ is defined by the
relation
\begin{equation*} 
  \frac{1}{\left\lvert *\sigma \right\rvert } \left\langle *\alpha, * \sigma
  \right\rangle = \kappa (\sigma ) \frac{1}{\left\lvert \sigma
    \right\rvert } \left\langle \alpha , \sigma \right\rangle,
\end{equation*} 
where $|\sigma |$ and $|*\sigma |$ are the volumes of these elements,
and $\kappa$ is the causality operator, which equals $ + 1 $ when
$\sigma $ is spacelike and $ - 1 $ otherwise.  (For more information
on alternative discrete Hodge operators, the reader may refer to,
e.g., \citealp{ArFaWi2006,AuKu2006,TaKeBo1999,WaWeToDeSc2006}.)

\subsubsection{Inner Product}
Define the \emph{$ \mathbb{L}^2 $ inner product} $ \left( \cdot ,
  \cdot \right) $ between two primal $k$-forms to be
\begin{align*} 
  \left( \alpha, \beta \right) &= \sum_{\sigma ^k} \kappa (\sigma)
  {n\choose k} \frac{ \left\lvert \operatorname{CH}(\sigma, *\sigma)
    \right\rvert }{ \left\lvert \sigma \right\rvert ^2 } \left\langle
    \alpha, \sigma \right\rangle \left\langle \beta, \sigma
  \right\rangle \\
  &= \sum_{\sigma ^k} \kappa (\sigma ) \frac{\left\lvert *\sigma
    \right\rvert }{\left\lvert \sigma \right\rvert } \left\langle
    \alpha, \sigma \right\rangle \left\langle \beta, \sigma
  \right\rangle
\end{align*} 
where the sum is taken over all $k$-dimensional elements $\sigma$, and
$\operatorname{CH}( \sigma , *\sigma )$ is the $n$-dimensional convex
hull of $ \sigma \cup *\sigma $ (see~\autoref{fig:primal-dual}).  The
final equality holds as a result of using the circumcentric dual,
since $\sigma$ and $ *\sigma$ are orthogonal to one another, and hence
$ \left\lvert \operatorname{CH}(\sigma , *\sigma ) \right\rvert =
{n\choose k}^{-1} \left\lvert \sigma \right\rvert \left\lvert *\sigma
\right\rvert $.  (Indeed, this is one of the advantages of using the
circumcentric dual, since one only needs to store volume information
about the primal and dual cells themselves, and not about these
primal-dual convex hulls.)  This inner product can be expressed in
terms of $ \alpha \wedge * \beta $, as in the continuous case, for a
particular choice of the discrete primal-dual wedge product;
see~\citet{DeHiMa2003}.

\begin{figure}
\includegraphics[width=\linewidth]{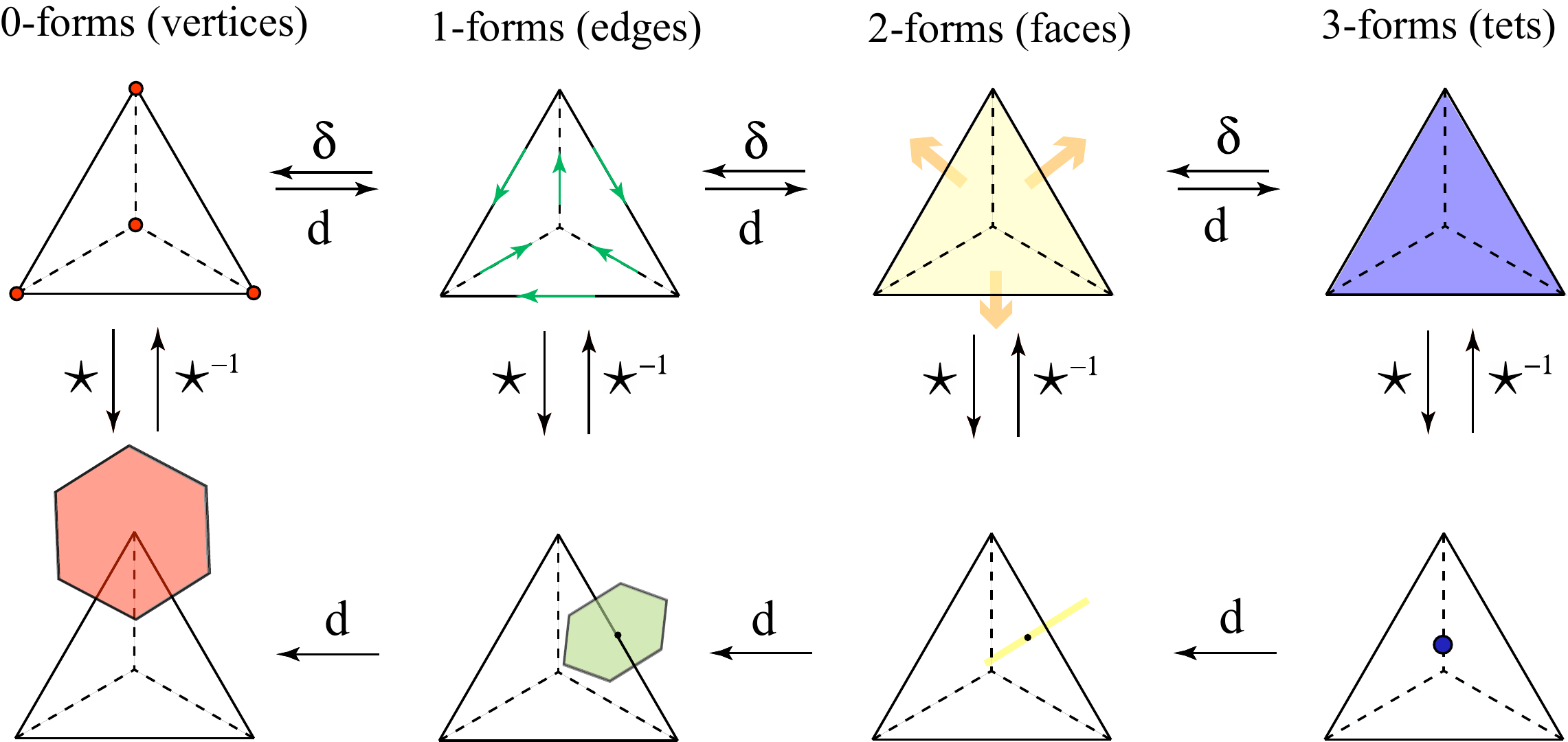}
\caption{This figure is an illustration of discrete differential forms
  and operators on a $3$-D simplicial mesh.  In the top row, we see
  how a discrete $k$-form lives on $k$-cells of the primal mesh, for $
  k = 0,1,2,3 $; the bottom row shows the location of the
  corresponding dual $(n-k)$-forms on the dual mesh.  The differential
  operators $ \mathrm{d} $ and $ \delta $ map ``horizontally'' between
  $k$ and $(k+1)$ forms, while the Hodge star $ * $ and its inverse $
  * ^{-1} $ map ``vertically'' between primal and dual forms.}
\label{fig:dualforms}
\end{figure}

Note that since we have already defined a discrete version of the
operators $ \mathrm{d} $ and $ * $, we immediately have a discrete
codifferential $ \delta $, with the same formal expression as given
previously.  See~\autoref{fig:dualforms} for a visual diagram of
primal and dual discrete forms, along with the corresponding operators
$ \mathrm{d}, *, \delta $, for the case where $K$ is a $3$-{D}
tetrahedral mesh.

\subsubsection{Implementing DEC}
DEC can be implemented simply and efficiently using linear algebra.  A
$k$-form $\alpha$ can be stored as a vector, where its entries are the
values of $\alpha$ on each $k$-cell of the mesh.  That is, given a
list of $k$-cells $ \sigma ^k _i $, the entries of the vector are $
\alpha _i = \left\langle \alpha , \sigma ^k _i \right\rangle $.  The
exterior derivative $ \mathrm{d} $, taking $k$-forms to $(k+1)$-forms,
is then represented as a matrix: in fact, it is precisely the
\emph{incidence matrix} between $k$-cells and $(k+1)$-cells in the
mesh, with sparse entries $\pm 1 $.  The Hodge star taking primal
$k$-forms to dual $(n-k)$-forms becomes a square matrix, and in the
case of the diagonal Hodge star, it is the diagonal matrix with
entries $ \kappa \left( \sigma ^k _i \right) \frac{\left\lvert
    *\sigma^k_i \right\rvert }{ \left\lvert \sigma^k_i \right\rvert}
$.  The discrete inner product is then simply the Hodge star matrix
taken as a quadratic form.

Because of this straightforward isomorphism between DEC and linear
algebra, problems posed in the language of DEC can take advantage of
existing numerical linear algebra codes.  For more details on
programming and implementation, refer to~\citet{ElSc2005}.

\subsection{Initial and Boundary Values with DEC}
\label{sec:boundaryDEC}
Particular care is required to properly enforce initial and boundary
conditions on the discrete spacetime boundary $ \partial K $.  For
example, in electromagnetism, we may wish to set initial conditions
for $E$ and $B$ at time $t_0$ --- but while $B$ is defined on
$\partial K$ at $t_0$, $E$ is not.  In fact, as we will see, $E$ lives
on edges that are extruded between the time slices $ t _0 $ and $ t _1
$, so unless we modify our definitions, we can only initialize $E$ at
the half-step $t_{ 1/2 }$.  (This half-step issue also arises with the
standard Yee scheme.)  There are some applications where it may be
acceptable to initialize $E$ and $B$ at separate times (for example,
when the fields are initialized randomly and integrated for a long
time to compute a resonance spectrum), but we wish to be able to
handle the more general case.  Although our previous exposition of DEC
thus far applies anywhere away from a boundary, notions as simple as
``dual cell'' need to be defined carefully on or near $ \partial K $.

\begin{figure}
\centerline{\includegraphics[width=0.7\linewidth]{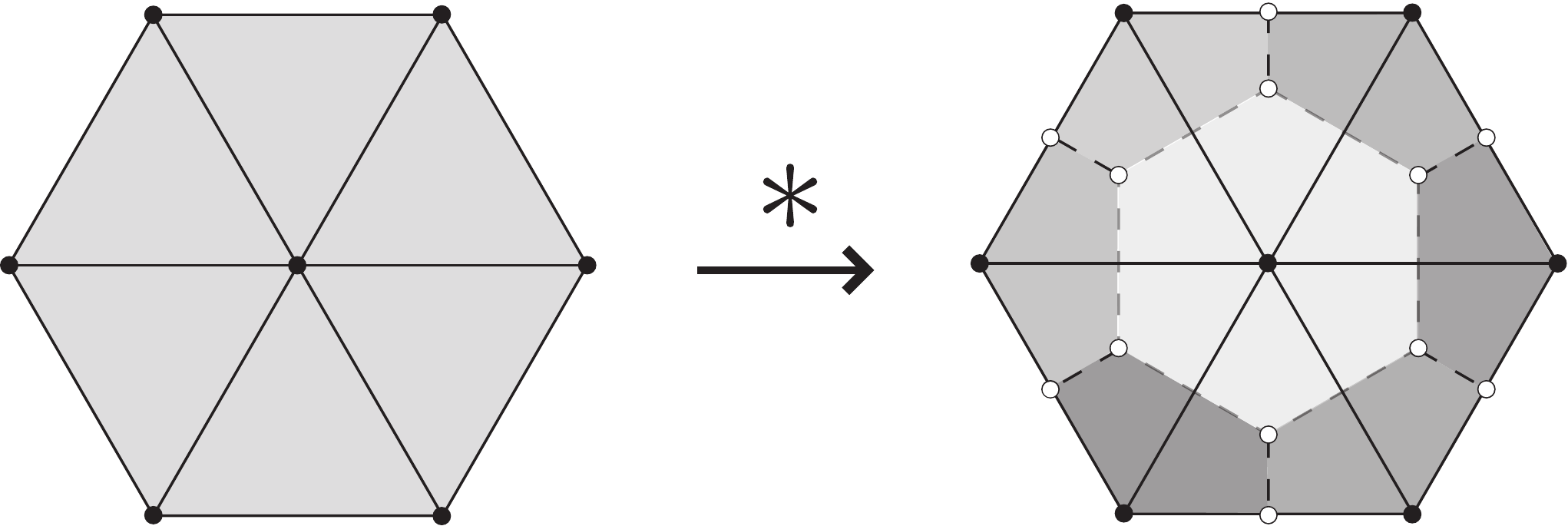}}
\caption[dual-boundary]{In this $2$-{D} example, the dual mesh is properly
  defined near the boundary by adding dual vertices on the boundary
  edges.  The restricted Voronoi cells of the primal boundary vertices
  (shaded at right) thus have boundaries containing both dual edges
  (dashed lines) \emph{and} primal boundary
  half-edges. } \label{fig:dual-boundary}
\end{figure}

For a primal mesh $K$, the dual mesh $*K$ is defined as \emph{the
  Voronoi dual of $K$ restricted to $K$}.  This truncates the portion
of the dual cells extending outside of $ K $; compare
\autoref{fig:dual-boundary} with the earlier
\autoref{fig:primal-dual}.  This new definition results in the
addition of a dual vertex at the circumcenter of each boundary
$(n-1)$-simplex, in addition to the interior $n$-simplices as
previously defined. To complete the dual mesh $*K$, we add a dual edge
between adjacent dual vertices on the boundary, as well as between
dual boundary vertices and their neighboring interior dual vertices,
and proceed similarly with higher-dimensional dual cells. For
intuition, one can imagine the $(n-1)$-dimensional boundary to be a
vanishingly thin $n$-dimensional shell.  That is, each boundary
$(k-1)$-simplex can be thought of as a prismal $k$-cell that has been
``squashed flat'' along the boundary normal direction.  This process
is quite similar to the use of ``ghost cells'' at the boundary, as is
commonly done for finite volume methods \citep[see][]{LeVeque2002}.
Note that these additional dual cells provide the boundary $\partial
K$ with its own dual mesh $* (\partial K)$. In fact, the boundary of
the dual is now equal to the dual of the boundary, i.e., $ \partial (*
K) = * (\partial K) $.
Returning to the example of initial conditions on $E$ and $B$, we
recall that $E$ is defined on extruded faces normal to the time slice
$ t _0 $. Therefore, thanks to the proper restriction of the Voronoi
diagram to the domain, we can now define $E$ on \emph{edges} in
$\partial K $ at time $ t _0 $, where these edges can be understood as
vanishingly thin faces (i.e., extruded between some $ t _{ - \epsilon
} $ and $ t _0 $ for $ \epsilon \rightarrow 0 $). Notice finally that
with this construction of $*K$, there is a dual relationship between
Dirichlet conditions on the dual mesh and Neumann conditions on the
primal mesh, e.g., between primal fields and dual fluxes, as expected.

\subsection{Discrete Integration by Parts with Boundary Terms}
With the dual mesh properly defined, dual forms can now be defined on
the boundary. Therefore, the discrete duality between $\mathrm{d} $
and $\delta$ can be generalized to include nonvanishing boundary
terms. If $\alpha$ is a primal $(k-1)$-form and $\beta$ is a primal
$k$-form, then
\begin{equation}
  \left(  \mathrm{d}  \alpha , \beta \right) 
  = \left(  \alpha , \delta \beta
  \right)  + \left\langle  \alpha \wedge *\beta , \partial K \right\rangle .
\label{discrete-ibp}
\end{equation}
In the boundary integral, $\alpha$ is still a primal $(k-1)$-form on $
\partial K $, while $ *\beta $ is an $(n-k)$-form taken on the
boundary dual $ *(\partial K) $.  Formula \eqref{discrete-ibp} is
readily proved using the familiar method of discrete ``summation by
parts,'' and thus agrees with the integration by parts formula for
smooth differential forms.

\section{Implementing Maxwell's Equations with DEC}
In this section, we explain how to obtain numerical algorithms for
solving Maxwell's equations with DEC.  To do so, we will proceed in
the following order.  First, we will find a sensible way to define the
discrete forms $F$, $G$, and $\mathcal{J}$ on a spacetime mesh.  Next,
we will use the DEC version of the operators $\mathrm{d}$ and $*$ to
obtain the discrete Maxwell's equations.  While we haven't yet shown
that these equations are variational in the discrete sense, we will
show later in \autoref{section:theory} that the Lagrangian derivation
of the smooth Maxwell's equations also holds with the DEC operators,
in precisely the same way.  Finally, we will discuss how these
equations can be used to define a numerical method for computational
electromagnetics.

In particular, for a rectangular grid, we will show that our setup
results in the traditional Yee scheme. For a general triangulation of
space with equal time steps, the resulting scheme will be Bossavit and
Kettunen's scheme.  We will then develop an AVI method, where each
spatial element can be assigned a different time step, and the time
integration of Maxwell's equations can be performed on the elements
asynchronously. Finally, we will comment on the equations for fully
generalized spacetime meshes, e.g., an arbitrary meshing of
$\mathbb{R}^{3,1} $ by $4$-simplices.

Note that the idea of discretizing Maxwell's equations using spacetime
cochains was mentioned in, e.g.,~\citet{Leok2004}, as well as in a
paper by~\citet{Wise2006} taking the more abstract perspective of
higher-level ``$p$-form'' versions of electromagnetism and category
theory.

\subsection{Rectangular Grid}
\label{section:dec-maxwell}

Suppose that we have a rectangular grid in $ \mathbb{R}^{3,1} $,
oriented along the axes $ (x, y, z, t) $.  To simplify this exposition
(although it is not necessary), let us also suppose that the grid has
uniform space and time steps $ \Delta x, \Delta y, \Delta z, \Delta t
$.  Note that the DEC setup applies directly to a non-simplicial
rectangular mesh, since an $n$-rectangle does in fact have a
circumcenter.

\subsubsection{Setup}
\begin{figure}
\centerline{\includegraphics[width=0.7\linewidth]{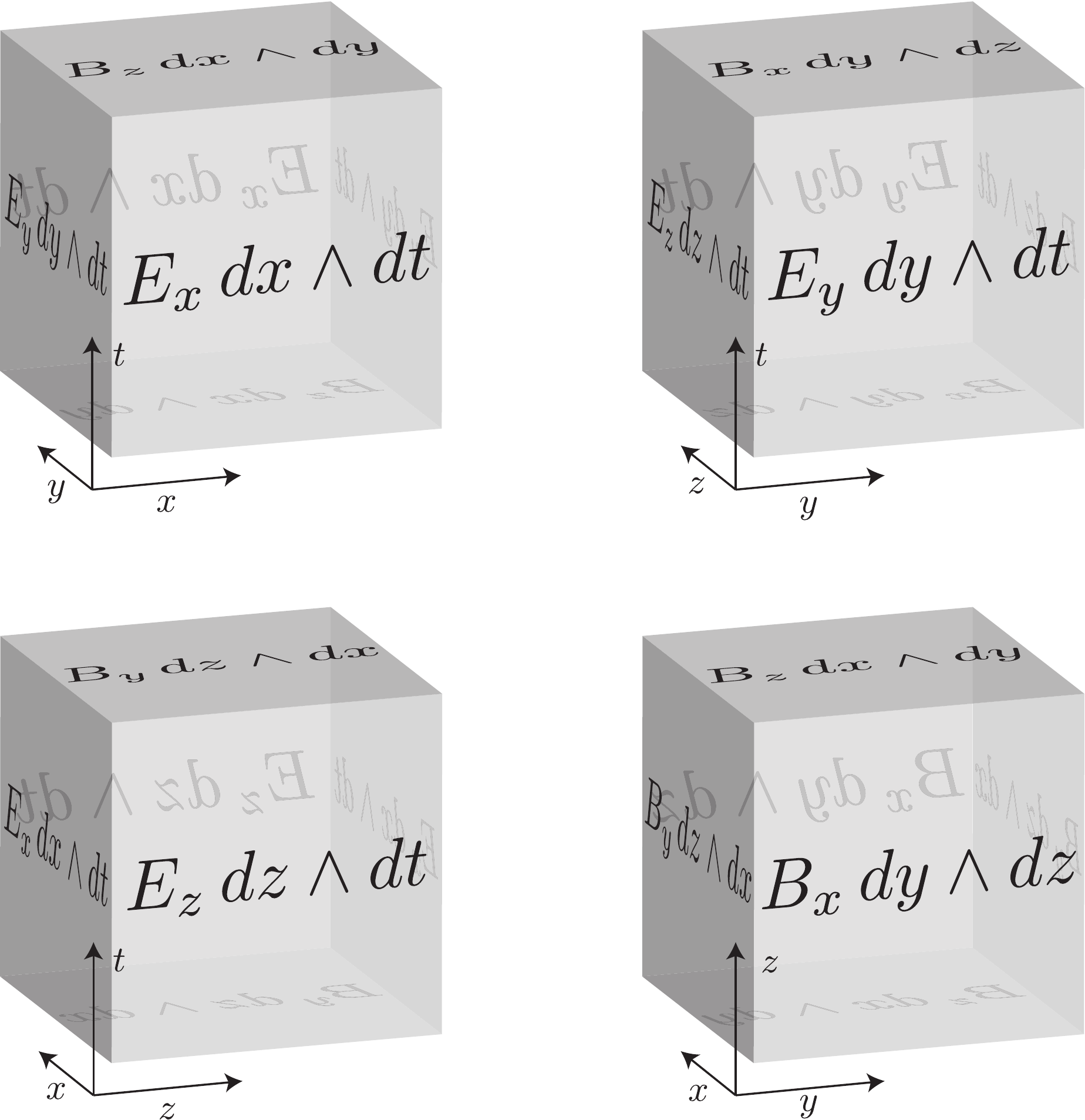}}
\caption[setupEB]{Values of $F$ are stored on the primal $2$-faces of
  a $4$-{D} rectangular grid.  Shown here are the three mixed
  space/time $3$-cells, and the one purely spatial $3$-cell (lower
  right).  } \label{fig:setupEB}
\end{figure}

Since $F$ is a $2$-form, its values should live on $2$-faces in this
grid. Following the continuous expression of $F$
\begin{align*} 
  F &= E _x \;\mathrm{d}x \wedge \mathrm{d} t + E _y \;\mathrm{d}y
  \wedge \mathrm{d} t + E _z \;\mathrm{d}z \wedge d t \\
  &\qquad + B _x \;\mathrm{d}y \wedge \mathrm{d} z + B _y \;\mathrm{d}z
  \wedge \mathrm{d} x + B _z \;\mathrm{d}x \wedge \mathrm{d} y,
\end{align*} 
and due to the tensor product nature of the regular grid, the exact
assignment of each $2$-face becomes simple: \emph{the six components
  of $F$ correspond precisely to the six types of $2$-faces in a
  $4$-{D} rectangular grid}.  Simply assign the values $ E _x \Delta x
\Delta t $ to faces parallel to the $ x t $-plane, $ E _y \Delta y
\Delta t $ to faces parallel to the $ y t $-plane, and $ E _z \Delta z
\Delta t $ to faces parallel to the $ z t $-plane.  Likewise, assign $
B _x \Delta y \Delta z $ to faces parallel to the $ yz $-plane, $ B _y
\Delta z \Delta x $ to faces parallel to the $ xz $-plane, and $ B _z
\Delta x \Delta y $ to faces parallel to the $ xy $-plane.  This is
pictured in \autoref{fig:setupEB}.

Let us look at these values on the faces of a typical $4$-rectangle $
\left[ x _k, x _{ k + 1 } \right] \times \left[ y _l , y _{ l + 1 }
\right] \times \left[ z _m , z _{ m + 1 } \right] \times \left[ t _n ,
  t _{ n + 1 } \right] $.  To simplify the notation, we can index each
value of $F$ by the midpoint of the $2$-face on which it lives: for
example, $ \left. F \right\rvert ^{n + \frac{1}{2} } _{ k +
  \frac{1}{2}, l, m} $ is stored on the face $ \left[ x _k, x _{ k + 1
  } \right] \times \left\{y _l\right\} \times \left\{z _m\right\}
\times \left[ t _n , t _{ n + 1 } \right] $, parallel to the $x
t$-plane.  Hence, the following values are assigned to the
corresponding faces:
\begin{align}
x t \text{-face}: &  \left. E _x \right\rvert ^{n + \frac{1}{2} } _{ k +
  \frac{1}{2}, l, m} \Delta x \Delta t \nonumber\\
yt \text{-face}: &  \left. E _y \right\rvert ^{n + \frac{1}{2} } _{k, l +
  \frac{1}{2}, m } \Delta y \Delta t \nonumber\\
zt \text{-face}: &  \left. E _z \right\rvert ^{n + \frac{1}{2} } _{k, l, m
+ \frac{1}{2} } \Delta z \Delta t \nonumber
\end{align}
\begin{align}
yz \text{-face}: &  \left. B _x \right\rvert ^n _{k, l + \frac{1}{2} , m
+ \frac{1}{2} } \Delta y \Delta z \nonumber\\
xz \text{-face}: &  \left. B _y \right\rvert ^n _{k + \frac{1}{2} , l, m
+ \frac{1}{2} } \Delta z \Delta x \nonumber\\
xy \text{-face}: &  \left. B _z \right\rvert ^n _{k + \frac{1}{2} , l +
  \frac{1}{2} , m} \Delta x \Delta y . \nonumber
\end{align}
We see that a ``staggered grid'' arises from the fact that $E$ and $B$
naturally live on $2$-faces, not at vertices or $4$-cells.

\subsubsection{Equations of Motion}
The discrete equations of motion are, as in the continuous case,
\begin{align*}
  \mathrm{d} F = 0, \qquad \mathrm{d} G = \mathcal{J} ,
\end{align*}
where now these equations are interpreted in the sense of DEC.  Let us
first look at the DEC interpretation of $ \mathrm{d} F $.  Since $
\mathrm{d} F $ is a discrete $3$-form, it takes values on the
$3$-faces of each $4$-rectangle.  Its values are as follows:
\begin{align}
xyt \text{-face}: &  - \left( \left. E _x \right\rvert ^{n + \frac{1}{2} } _{ k +
  \frac{1}{2}, l + 1, m} - \left. E _x \right\rvert ^{n + \frac{1}{2} } _{ k +
  \frac{1}{2}, l, m} \right) \Delta x \Delta t \nonumber \\
&  + \left( \left. E _y \right\rvert ^{n + \frac{1}{2} } _{k + 1, l +
  \frac{1}{2}, m } - \left. E _y \right\rvert ^{n + \frac{1}{2} } _{k, l +
  \frac{1}{2}, m } \right) \Delta y \Delta t  \nonumber \\
&  + \left( \left. B _z \right\rvert ^{n+1} _{k + \frac{1}{2} , l +
  \frac{1}{2} , m} - \left. B _z \right\rvert ^n _{k + \frac{1}{2} , l +
  \frac{1}{2} , m} \right) \Delta x \Delta y \nonumber
  \end{align}
  \begin{align}
xzt \text{-face}:
&  - \left( \left. E _x \right\rvert ^{n + \frac{1}{2} } _{ k +
  \frac{1}{2}, l, m+1} - \left. E _x \right\rvert ^{n + \frac{1}{2} } _{ k +
  \frac{1}{2}, l, m} \right) \Delta x \Delta t \nonumber \\
&  + \left( \left. E _z \right\rvert ^{n + \frac{1}{2} } _{k+1, l, m
+ \frac{1}{2} } - \left. E _z \right\rvert ^{n + \frac{1}{2} } _{k, l, m
+ \frac{1}{2} } \right) \Delta z \Delta t \nonumber \\
&  - \left( \left. B _y \right\rvert ^{n+1} _{k + \frac{1}{2} , l, m
+ \frac{1}{2} } - \left. B _y \right\rvert ^n _{k + \frac{1}{2} , l, m
+ \frac{1}{2} } \right) \Delta x \Delta z \nonumber
\end{align}
\begin{align}
yzt \text{-face}:
&  - \left( \left. E _y \right\rvert ^{n + \frac{1}{2} } _{k, l +
  \frac{1}{2}, m +1} - \left. E _y \right\rvert ^{n + \frac{1}{2} } _{k, l +
  \frac{1}{2}, m } \right) \Delta y \Delta t \nonumber \\
&  + \left( \left. E _z \right\rvert ^{n + \frac{1}{2} } _{k, l+1, m
+ \frac{1}{2} } - \left. E _z \right\rvert ^{n + \frac{1}{2} } _{k, l, m
+ \frac{1}{2} } \right) \Delta z \Delta t \nonumber \\
&  + \left( \left. B _x \right\rvert ^{n+1} _{k, l + \frac{1}{2} , m
+ \frac{1}{2} } - \left. B _x \right\rvert ^n _{k, l + \frac{1}{2} , m
+ \frac{1}{2} } \right) \Delta y \Delta z \nonumber
\end{align}
\begin{align}
xyz \text{-face}:
&  \left( \left. B _x \right\rvert ^n _{k+1, l + \frac{1}{2} , m
+ \frac{1}{2} } - \left. B _x \right\rvert ^n _{k, l + \frac{1}{2} , m
+ \frac{1}{2} } \right) \Delta y \Delta z \nonumber \\
&  + \left( \left. B _y \right\rvert ^n _{k + \frac{1}{2} , l+1, m
+ \frac{1}{2} } - \left. B _y \right\rvert ^n _{k + \frac{1}{2} , l, m
+ \frac{1}{2} } \right) \Delta x \Delta z \nonumber \\
&  + \left( \left. B _z \right\rvert ^n _{k + \frac{1}{2} , l +
  \frac{1}{2} , m+1} - \left. B _z \right\rvert ^n _{k + \frac{1}{2} , l +
  \frac{1}{2} , m} \right) \Delta x \Delta y \nonumber
\end{align}
Setting each of these equal to zero, we arrive at the following
four equations:
\begin{align}
& \frac{ \left. B _x \right| ^{n+1} _{k, l + \frac{1}{2} , m
+ \frac{1}{2} } - \left. B _x \right| ^n _{k, l + \frac{1}{2} , m
+ \frac{1}{2} } }{ \Delta t } = \nonumber \\
& \qquad
\frac{ \left. E _y \right| ^{n + \frac{1}{2} } _{k, l +
  \frac{1}{2}, m +1} - \left. E _y \right| ^{n + \frac{1}{2} } _{k, l +
  \frac{1}{2}, m } }{ \Delta z } - \frac{ \left. E _z \right| ^{n + \frac{1}{2} } _{k, l+1, m
+ \frac{1}{2} } - \left. E _z \right| ^{n + \frac{1}{2} } _{k, l, m
+ \frac{1}{2} } }{ \Delta y } \nonumber
\end{align}
\begin{align}
&  \frac{ \left. B _y \right| ^{n+1} _{k + \frac{1}{2} , l, m
+ \frac{1}{2} } - \left. B _y \right| ^n _{k + \frac{1}{2} , l, m
+ \frac{1}{2} } }{ \Delta t } = \nonumber \\
& \qquad
\frac{ \left. E _z \right| ^{n + \frac{1}{2} } _{k+1, l, m
+ \frac{1}{2} } - \left. E _z \right| ^{n + \frac{1}{2} } _{k, l, m
+ \frac{1}{2} } }{ \Delta x } - \frac{ \left. E _x \right| ^{n + \frac{1}{2} } _{ k +
  \frac{1}{2}, l, m+1} - \left. E _x \right| ^{n + \frac{1}{2} } _{ k +
  \frac{1}{2}, l, m} }{ \Delta z } \nonumber
  \end{align}
  \begin{align}
&  \frac{ \left. B _z \right| ^{n+1} _{k + \frac{1}{2} , l +
  \frac{1}{2} , m} - \left. B _z \right| ^n _{k + \frac{1}{2} , l +
  \frac{1}{2} , m} }{ \Delta t } =  \nonumber \\
& \qquad
\frac{ \left. E _x \right| ^{n + \frac{1}{2} } _{ k +
  \frac{1}{2}, l + 1, m} - \left. E _x \right| ^{n + \frac{1}{2} } _{ k +
  \frac{1}{2}, l, m} }{ \Delta y } - \frac{ \left. E _y \right| ^{n + \frac{1}{2} } _{k + 1, l +
  \frac{1}{2}, m } - \left. E _y \right| ^{n + \frac{1}{2} } _{k, l +
  \frac{1}{2}, m } }{ \Delta x } \nonumber
  \end{align}
  and
  \begin{equation}
    \begin{split}
    & \frac{ \left. B _x \right| ^n _{k+1, l + \frac{1}{2} , m +
          \frac{1}{2} } - \left. B _x \right| ^n _{k, l + \frac{1}{2}
          , m + \frac{1}{2} } }{ \Delta x } + \frac{ \left. B _y
        \right| ^n _{k + \frac{1}{2} , l+1, m + \frac{1}{2} } -
        \left. B _y \right| ^n _{k + \frac{1}{2} , l, m + \frac{1}{2}
      } }{ \Delta y } \\ & \qquad + \frac{ \left. B _z \right| ^n _{k
          + \frac{1}{2} , l + \frac{1}{2} , m+1} - \left. B _z \right|
        ^n _{k + \frac{1}{2} , l + \frac{1}{2} , m} }{ \Delta z } = 0 .
    \end{split}
    \label{yeedivb}
  \end{equation}
These equations are the discrete version of the equations
\begin{equation*} 
  \partial _t \mathbf{B}  = - \nabla \times \mathbf{E} , \qquad 
  \nabla \cdot \mathbf{B}  = 0.
\end{equation*} 
Moreover, since $E$ and $B$ are differential forms, this can also be
seen as a discretization of the {\em integral version} of Maxwell's
equations as well!  Because DEC satisfies a discrete Stokes' theorem,
this automatically preserves the equivalence between the differential
and integral formulations of electromagnetism.

\begin{figure}
\centerline{\includegraphics[width=0.7\linewidth]{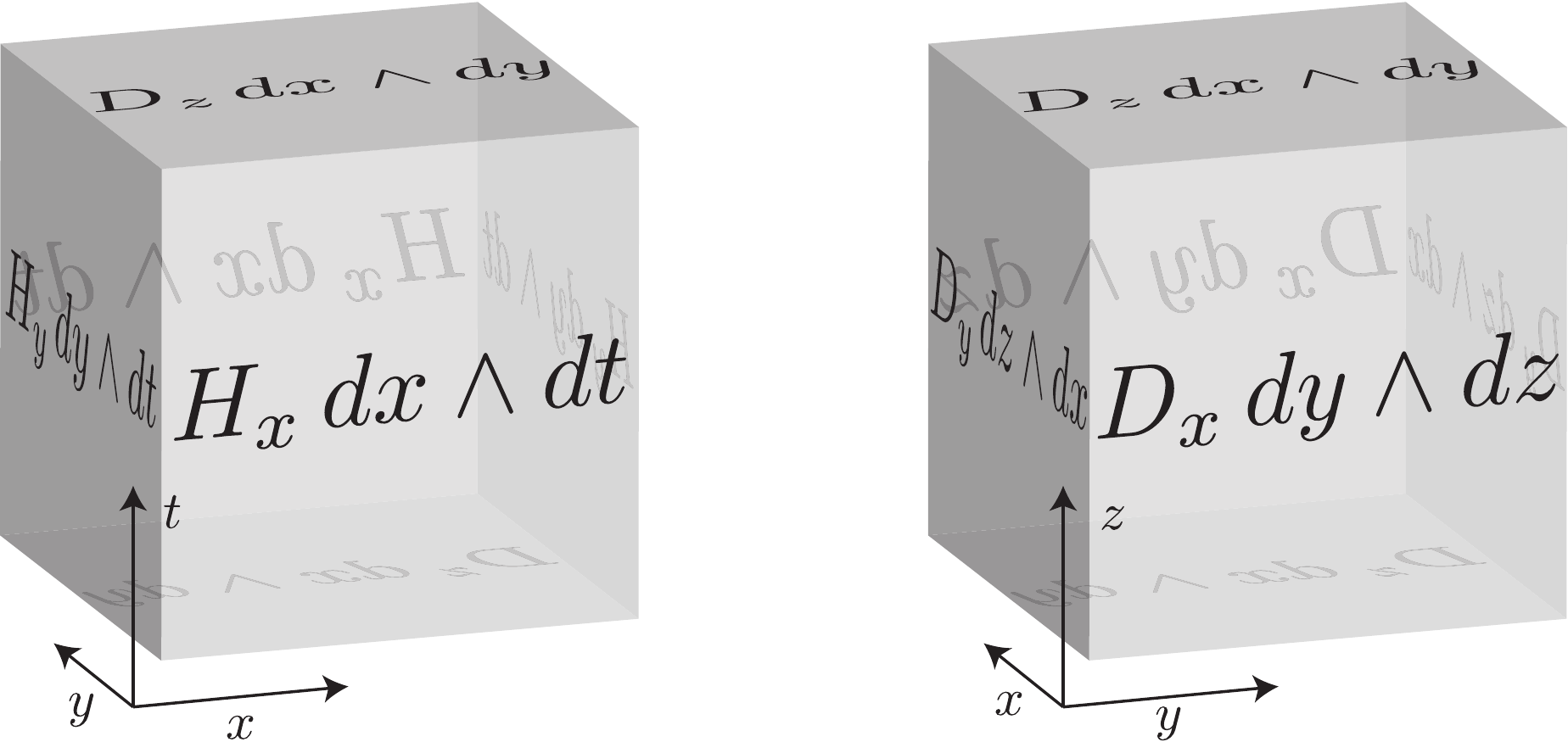}}
\caption[setupHD]{Values of $G = *F$ are stored on dual $2$-faces in a
  rectangular grid.  Shown here are a mixed space/time dual $3$-cell
  (left), corresponding to a spacelike primal edge; and a purely
  spatial dual $3$-cell (right), corresponding to a timelike primal
  edge.  There are also two other mixed space/time cells, as in
  \autoref{fig:setupEB}, that are not shown here.} \label{fig:setupHD}
\end{figure}

Doing the same with the equation $ \mathrm{d} G = \mathcal{J} $,
evaluating on dual $3$-faces this time, we arrive at four more
equations:
\begin{align}
& \frac{ \left. D _x \right\rvert ^{n + \frac{1}{2} } _{ k +
  \frac{1}{2}, l, m} - \left. D _x \right\rvert ^{n - \frac{1}{2} } _{ k +
  \frac{1}{2}, l, m} }{ \Delta t } =\nonumber \\
& \qquad \frac{ \left. H _z \right\rvert ^n _{k + \frac{1}{2} , l + \frac{1}{2} ,
    m} - \left. H _z \right\rvert ^n _{k + \frac{1}{2} , l - \frac{1}{2} ,
    m} }{ \Delta y } - \frac{ \left. H _y \right\rvert ^n _{k + \frac{1}{2}
    , l, m + \frac{1}{2} } - \left. H _y \right\rvert ^n _{k + \frac{1}{2}
    , l, m - \frac{1}{2} } }{ \Delta z } - \left. J _x \right\rvert ^n
_{ k + \frac{1}{2}, l , m } \nonumber
    \end{align}
    \begin{align}
      & \frac{ \left. D _y \right\rvert ^{n + \frac{1}{2} } _{k, l +
          \frac{1}{2}, m } - \left. D _y \right\rvert ^{n - \frac{1}{2} }
        _{k,
          l + \frac{1}{2}, m } }{ \Delta t } = \nonumber \\
      & \qquad \frac{ \left. H _x \right\rvert ^n _{k, l + \frac{1}{2} , m
          + \frac{1}{2} } - \left. H _x \right\rvert ^n _{k, l +
          \frac{1}{2} , m - \frac{1}{2} } }{ \Delta z } - \frac{
        \left. H _z \right\rvert ^n _{k + \frac{1}{2} , l + \frac{1}{2} ,
          m} - \left. H _z \right\rvert ^n _{k - \frac{1}{2} , l +
          \frac{1}{2} , m} }{ \Delta x } - \left. J _y \right\rvert ^n
      _{ k , l + \frac{1}{2} , m } \nonumber
    \end{align}
    \begin{align}
      & \frac{ \left. D _z \right\rvert ^{n + \frac{1}{2} } _{k, l, m +
          \frac{1}{2} } - \left. D _z \right\rvert ^{n - \frac{1}{2} } _{k,
          l, m +
          \frac{1}{2} }}{ \Delta t } = \nonumber \\
      & \qquad \frac{ \left. H _y \right\rvert ^n _{k + \frac{1}{2} , l, m
          + \frac{1}{2} } - \left. H _y \right\rvert ^n _{k - \frac{1}{2} ,
          l, m + \frac{1}{2} } }{ \Delta x } - \frac{ \left. H _x
        \right\rvert ^n _{k, l + \frac{1}{2} , m + \frac{1}{2} } - \left. H
          _x \right\rvert ^n _{k, l - \frac{1}{2} , m + \frac{1}{2} } }{
        \Delta y } - \left. J _z \right\rvert ^n _{ k , l , m +
        \frac{1}{2} } \nonumber
    \end{align}
    and
    \begin{equation} 
      \begin{split}
        & \frac{ \left. D _x \right\rvert ^{n + \frac{1}{2} } _{ k +
            \frac{1}{2}, l, m} - \left. D _x \right\rvert ^{n +
            \frac{1}{2} } _{ k - \frac{1}{2}, l, m} }{ \Delta x } +
        \frac{\left. D _y \right\rvert ^{n + \frac{1}{2} } _{k, l +
            \frac{1}{2}, m } - \left. D _y \right\rvert ^{n +
            \frac{1}{2} } _{k, l - \frac{1}{2}, m } }{ \Delta y
        }  \\
        & \qquad + \frac{ \left. D _z \right\rvert ^{n +
            \frac{1}{2} } _{k, l, m + \frac{1}{2} } - \left. D _z
          \right\rvert ^{n + \frac{1}{2} } _{k, l, m - \frac{1}{2} }
        }{ \Delta z } = \left. \rho \right\rvert ^{ n +\frac{1}{2} }
        _{ k, l, m } .
      \end{split}
      \label{yeedivd}
    \end{equation}
This results from storing $ G $ on the dual grid, as shown in
\autoref{fig:setupHD}.  This set of equations is the discrete version
of
\begin{equation*} 
  \partial _t \mathbf{D}   =    \nabla \times  \mathbf{H}  -
  \mathbf{J}  , \qquad 
  \nabla \cdot \mathbf{D}  =   \rho .
\end{equation*} 

After eliminating the redundant divergence equations \eqref{yeedivb}
and \eqref{yeedivd} (see \autoref{section:gauge} for details) and
making the substitutions $ \mathbf{D} = \epsilon \mathbf{E} $, $
\mathbf{B} = \mu \mathbf{H} $, the remaining equations are precisely
the Yee scheme, as formulated in~\citet[pp. 67--68]{BoRyIn2005}.

\subsection{Unstructured Spatial Mesh with Uniform Time Steps}
\label{Unstructured}

We now consider the case of an unstructured grid in space, but with
uniform steps in time as advocated in, e.g., \citet{BoKe1999}.
Suppose that, instead of a rectangular grid for both space and time,
we have an arbitrary space discretization on which we would like to
take uniform time steps.  (For example, we may be given a tetrahedral
mesh of the spatial domain.)  This mesh contains two distinct types of
$2$-faces.  First, there are triangular faces that live entirely in
the space mesh at a single position in time.  Every edge of such a
face is spacelike---that is, it has positive length---so the causality
operator defined in \autoref{sec:dec} takes the value $\kappa = 1 $.
Second, there are rectangular faces that live between time steps.
These faces consist of a single spacelike edge extruded by one time
step. Because they have one timelike edge, these faces satisfy $
\kappa = -1 $.  Again, the circumcentric-dual DEC framework applies
directly to this type of mesh, since the prismal extrusion of a
$3$-simplex still has a circumcenter.

\subsubsection{Setup}
\begin{figure}
\centerline{\includegraphics[width=0.7\linewidth]{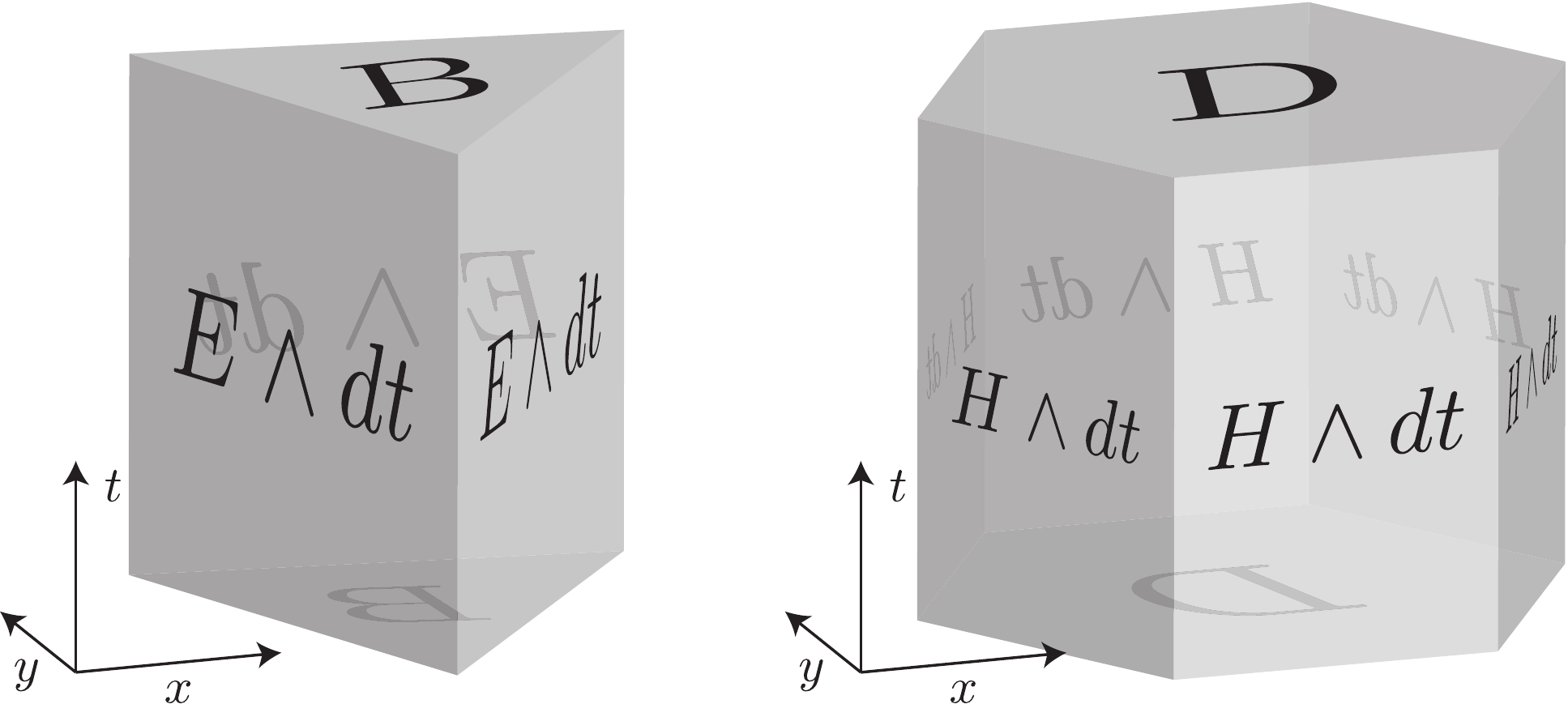}}
\caption[Setup]{For an unstructured spatial mesh, $F$ is stored on
  primal $2$-faces (left), while $ G = *F $ is stored on dual
  $2$-faces (right).  Shown here are the values on mixed space/time
  $3$-cells.  (The purely spatial $3$-cells, which correspond to the
  divergence equations and do not contribute to the equations of
  motion, are not shown.) } \label{fig:setupUnstruct}
\end{figure}

Again, we can characterize the discrete values of $F$ by looking at
the continuous expression
\begin{equation}
  F = E \wedge \mathrm{d}  t + B. \nonumber
\end{equation}
Therefore, let us assign $B$ to the purely spacelike faces and $E
\Delta t$ to the mixed space/time faces.  Looking at $ G = *F $ shows
that mixed dual faces should store $H \Delta t$ and spacelike dual
faces should store $D$; see \autoref{fig:setupUnstruct}.

\subsubsection{Equations of Motion}
As in~\citet{Bossavit1998}, we can store the values of each
differential form over every spatial element in an array, using the
method described in \autoref{sec:dec}.  This leads to the arrays $ B^n
$ and $ H^n $ at whole time steps $n$, and $ E^{n+1/2} $ and $ {D}
^{n+1/2} $ at half time steps.  Let $\mathrm{d} _1 $ denote the
edges-to-faces incidence matrix for the spatial domain.  That is, $
\mathrm{d} _1 $ is the matrix corresponding to the discrete exterior
derivative, taken only in space, from primal $1$-forms to primal
$2$-forms.  Similarly, the transpose $ \mathrm{d} _1 ^T $ corresponds
to the exterior derivative from spatial dual $1$-forms to dual
$2$-forms. Then the equation $ \mathrm{d} F = 0 $, evaluated on all
prismal $3$-faces becomes
\begin{equation}
  \frac{ B^{n+1} - B ^n }{ \Delta t } = - \mathrm{d} _1  E^{n + 1/2 }. \nonumber
\end{equation}
Likewise, the equation $ \mathrm{d} G = \mathcal{J} $, evaluated on
all space/time $3$-faces in the dual mesh, becomes
\begin{equation}
  \frac{ D ^{ n + 1/2 } - D ^{ n - 1/2 } }{ \Delta t } =
  \mathrm{d} _1 ^T  H ^n - J ^n . \nonumber
\end{equation}
We can also evaluate $ \mathrm{d} F = 0 $ and $ \mathrm{d} G =
\mathcal{J} $ on spacelike $3$-faces, e.g, tetrahedra; these simply
yield the discrete versions of the divergence conditions for $B$ and
$D$, which can be eliminated.

Therefore, the DEC scheme for such a mesh is equivalent to Bossavit
and Kettunen's Yee-like scheme; additionally, when the spatial mesh
is taken to be rectangular, this integrator reduces to the standard
Yee scheme. However, we now have solid foundations to extend this
integrator to handle asynchronous updates for improved efficiency.

\subsection{Unstructured Spatial Mesh with Asynchronous Time Steps}
\label{AVI}
Instead of picking the same time step size for every element of the
spatial mesh, as in the previous two sections, it is often more
efficient to assign each element its own, optimized time step, as done
in~\citet{LeMaOrWe2003} for problems in elastodynamics. In this case,
rather than the entire mesh evolving forward in time simultaneously,
individual elements advance one-by-one, asynchronously---hence the
name \emph{asynchronous variational integrator} (AVI). As we will
prove in \autoref{section:theory}, this asynchronous update process
will maintain the variational nature of the integration scheme. Here,
we again allow the spatial mesh to be unstructured.

\subsubsection{Setup}
After choosing a primal space mesh, assign each spatial $2$-face
(e.g., triangle) $\sigma$ its own discrete time set
\begin{equation}
  \Theta _\sigma = \left\{ t ^0 _\sigma < \cdots < t ^{ N _\sigma }
    _\sigma 
  \right\}. \nonumber
\end{equation}
For example, one might assign each face a fixed time step size $
\Delta t _\sigma = t ^{n+1} _\sigma - t ^n _\sigma $, taking equal
time steps \emph{within} each element, but with $ \Delta t $ varying
\emph{across} elements.  We further require for simplicity of
explanation that, except for the initial time, no two faces take the
same time step: that is, $ \Theta _\sigma \cap \Theta _{\sigma
  ^\prime} = \{t_0\} $ for $ \sigma \neq \sigma ^\prime $.

\begin{figure}
\centerline{\includegraphics[width=0.7\linewidth]{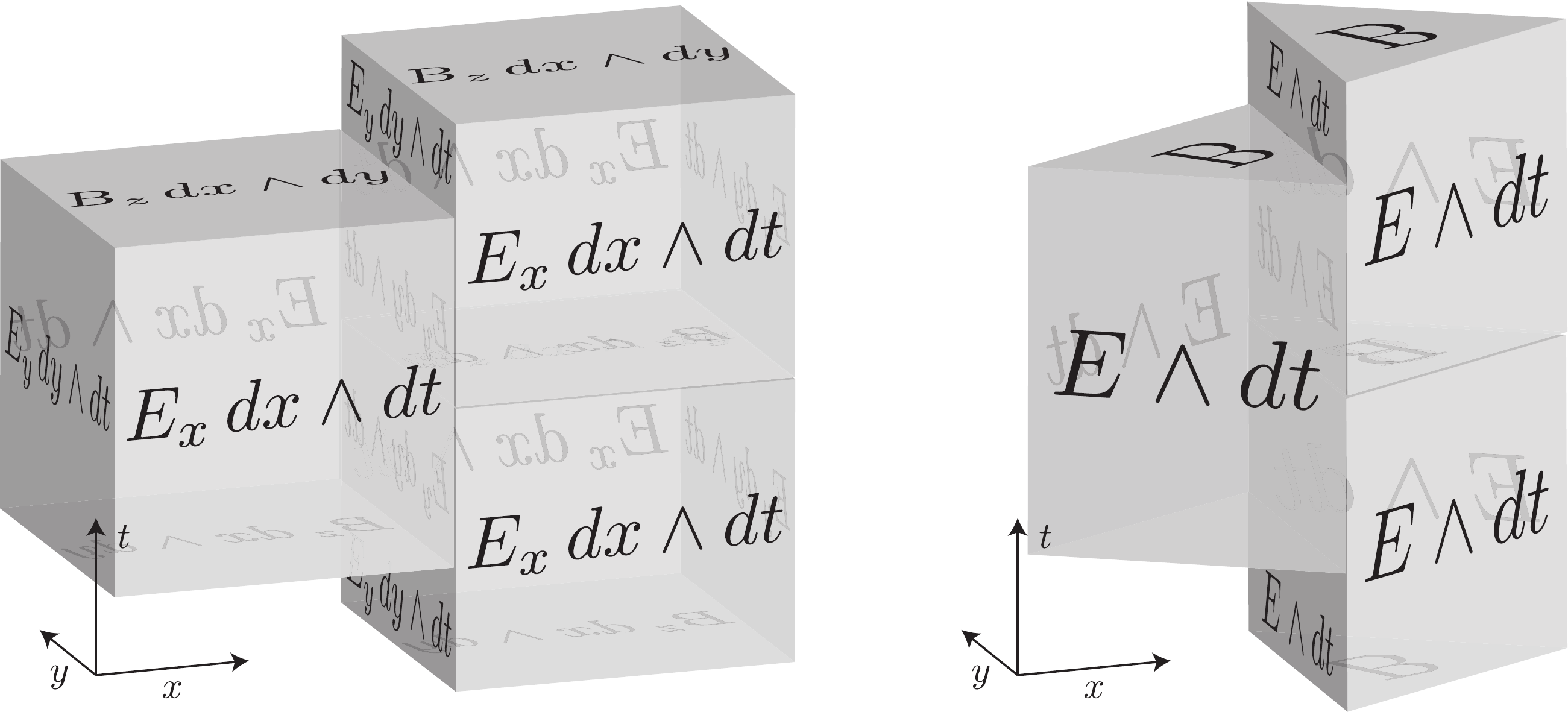}}
\caption[setupAVI]{Shown here is part of an AVI mesh, for a
  rectangular spatial mesh (left) and for an unstructured spatial mesh
  (right).  The different heights of the spacetime prisms reflect the
  fact that elements can take different time steps from one another.
  Moreover, these time steps can be asynchronous, as seen in the
  mismatch between the horizontal faces.} \label{fig:setupAVI}
\end{figure}

In order to keep proper time at the edges $e$ where multiple faces with
different time sets meet, we let
\begin{equation}
  \Theta _e = \bigcup _{ \sigma \ni e } \Theta _\sigma = \left\{ t ^0 _e \leq
    \cdots \leq t ^{ N _e } _e \right\}. \nonumber 
\end{equation}
Therefore the mixed space-time $2$-faces, which correspond to the edge
$e$ extruded over a time step, are assigned the set of intermediate
times
\begin{equation}
\Theta _{e} ^\prime = \left\{ t ^{1/2} _e \leq \cdots \leq t ^{ N _e
- 1/2 } _e \right\}, \nonumber 
\end{equation}
where $ t ^{k + 1/2} _e = (t ^{k + 1} _e + t ^{k} _e)/2 $.  The values
stored on a primal AVI mesh are shown in \autoref{fig:setupAVI}.

Since $ \Theta _e \supset \Theta _\sigma $ when $ e \subset \sigma $,
each spatial edge $e$ takes more time steps than any one of its
incident faces $\sigma$; as a result, it is not possible in general to
construct a circumcentric dual on the entire spacetime AVI mesh, since
the mesh is not prismal and hence the circumcenter may not exist.
Instead, we find the circumcentric dual to the \emph{spatial} mesh,
and assign same time steps to the primal and dual elements
\begin{equation*}
  \Theta _{*\sigma} = \Theta _\sigma , \qquad 
  \Theta _{*e} = \Theta _{e}.
\end{equation*} 
This results in well-defined primal and dual cells for each
$2$-element in spacetime, and hence a Hodge star for this order.  (A
Hodge star on forms of different order is not needed to formulate
Maxwell's equations.)

\subsubsection{Equations of Motion}
The equation $ \mathrm{d} F = 0 $, evaluated on a mixed space/time
$3$-cell, becomes
\begin{equation}
  \label{avi-update-b} \frac{ B^{n+1}_\sigma  - B ^n_\sigma  }{
    t^{n+1}_\sigma  - t^n_\sigma } = - \mathrm{d} _1  \sum \left\{
    E^{m+1/2}_e : t^n_\sigma  < t^{m+1/2}_e < t^{n+1}_\sigma \right\} .
\end{equation}
Similarly, the equation $ \mathrm{d} G = \mathcal{J} $ becomes
\begin{equation}
\label{avi-update-e}
\frac{ D ^{ m + 1/2 }_e - D ^{ m - 1/2 }_e }{ t^{m+1/2}_e -
  t^{m-1/2}_e } = \mathrm{d} _1 ^T  \left( H ^n_\sigma  \,
  \mathbb{1}_{ \left\{t^n_\sigma = t^m_e \right\}} \right) - J ^m _e ,
\end{equation}
where $ \mathbb{1}_{ \left\{t^n_\sigma = t^m_e \right\}} $ equals $1$
when face $\sigma$ has $ t^n_\sigma = t^m_e $ for some $n$, and $0$
otherwise.  (That is, the indicator function ``picks out'' the
incident face that lives at the same time step as this edge.)

Solving an initial value problem can then be summarized by the
following update loop:
\begin{enumerate}
\item Pick the minimum time $ t^{n+1}_\sigma $ where $ B^{n+1}_\sigma $ has not
  yet been computed.
\item Advance $B^{n+1}_\sigma $ according to \autoref{avi-update-b}.
\item Update $ H ^{ n + 1 } _\sigma = *_\mu ^{-1} B ^{n+1}_\sigma $.
\item Advance $D ^{ m + 3/2 }_e$ on neighboring edges $ e \subset
  \sigma $ according to \autoref{avi-update-e}.
\item Update $E^{m+3/2}_e = *_\epsilon ^{-1} D ^{ m + 3/2 }_e $.
\end{enumerate}

\subsubsection{Iterative Time Stepping Scheme}
As detailed in~\citet{LeMaOrWe2003} for elastodynamics, the explicit
AVI update scheme can be implemented by selecting mesh elements from a
priority queue, sorted by time, and iterating forward.  However, as
written above, the scheme is not strictly iterative, since
\autoref{avi-update-e} depends on past values of $E$.  This can be
easily fixed by rewriting the AVI scheme to advance in the variables
$A$ and $E$ instead, where the potential $A$ effectively stores the
cumulative contribution of $ E $ to the value of $B $ on neighboring
faces. Compared to the AVI for elasticity, $A$ plays the role of the
positions $\mathbf{x}$, while $E$ plays the role of the (negative)
velocities $\dot{\mathbf{x}} $.  The algorithm is given as pseudocode
in \autoref{fig:pseudocode}. Note that if all elements take uniform
time steps, the AVI reduces to the Bossavit--Kettunen scheme.

\begin{figure}
\centering
\begin{minipage}[c]{0.9\columnwidth} %
\hrule \vspace*{2mm}%
\begin{algorithmic}
\STATE \COMMENT{\textsc{Initialize fields and priority queue}}%
\FOR{each spatial edge $e$}%
    \STATE $A_e \gets A_e^0, \; E_e \gets E_e ^{1/2}, \; \tau_e \gets
    t_0$ \emph{// Store initial field values and times}
\ENDFOR
\FOR{each spatial face $\sigma$}%
    \STATE $\tau_\sigma \gets t_0$
    \STATE Compute the next update time $t_\sigma^1$
    \STATE $Q$.push$(t_\sigma^1,\sigma)$ \emph{// Push element onto queue with
      its next update time}
\ENDFOR
\vspace*{3mm}
\COMMENT{\textsc{Iterate forward in time until the priority queue is empty}}%
\REPEAT
    \STATE $(t,\sigma) \gets$ $Q$.pop() \emph{// Pop next element $\sigma$ and time $t$ from queue}
    \FOR{each edge $e$ of element $\sigma$}%
        \STATE $A_e \gets A_e - E_e(t-\tau_e)$ \emph{// Update
          neighboring values of $A$ at time $t$}
    \ENDFOR
    \IF{$t<\mbox{final-time}$} \STATE $ B _\sigma \gets \mathrm{d}_1 A
    _e $ \STATE $ H _\sigma \gets *_\mu B _\sigma$ \STATE $ D _e \gets
    *_\epsilon E _e $ \STATE $ D _e \gets D _e + \mathrm{d}_1(e,
    \sigma) H _\sigma (t - \tau _\sigma )$ \STATE $ E _e \gets
    *_\epsilon D _e $ \STATE $ \tau _\sigma \gets t $ \emph{// Update
      element's time} \STATE Compute the next update time $ t _\sigma
    ^{\text{next}} $ \STATE $Q$.push$(t _\sigma ^{\text{next}}, \sigma
    )$ \emph{// Schedule $\sigma$ for next update}
    \ENDIF
\UNTIL{($Q$.isEmpty())}
\end{algorithmic}\vspace*{2mm}
\hrule
\end{minipage}\vspace*{-2.5mm}
\caption{Pseudocode for our Asynchronous Variational Integrator,
  implemented using a priority queue data structure for storing and
  selecting the elements to be
  updated. \label{fig:pseudocode}}\vspace*{-5mm}
\end{figure}

\subsubsection{Numerical Experiments}
We first present a simple numerical example demonstrating the good
energy behavior of our asynchronous integrator.  The AVI was used to
integrate in time over a $2$-{D} rectangular cavity with perfectly
electrically conducting (PEC) boundaries, so that $E$ vanishes at the
boundary of the domain. $E$ was given random values at the initial
time, so as to excite all frequency modes, and integrated for 8
seconds.  Each spatial element was given a time step equal to 1/10 of
the stability-limiting time step determined by the CFL condition.

\begin{figure}
\centerline{\includegraphics[width=0.7\linewidth]{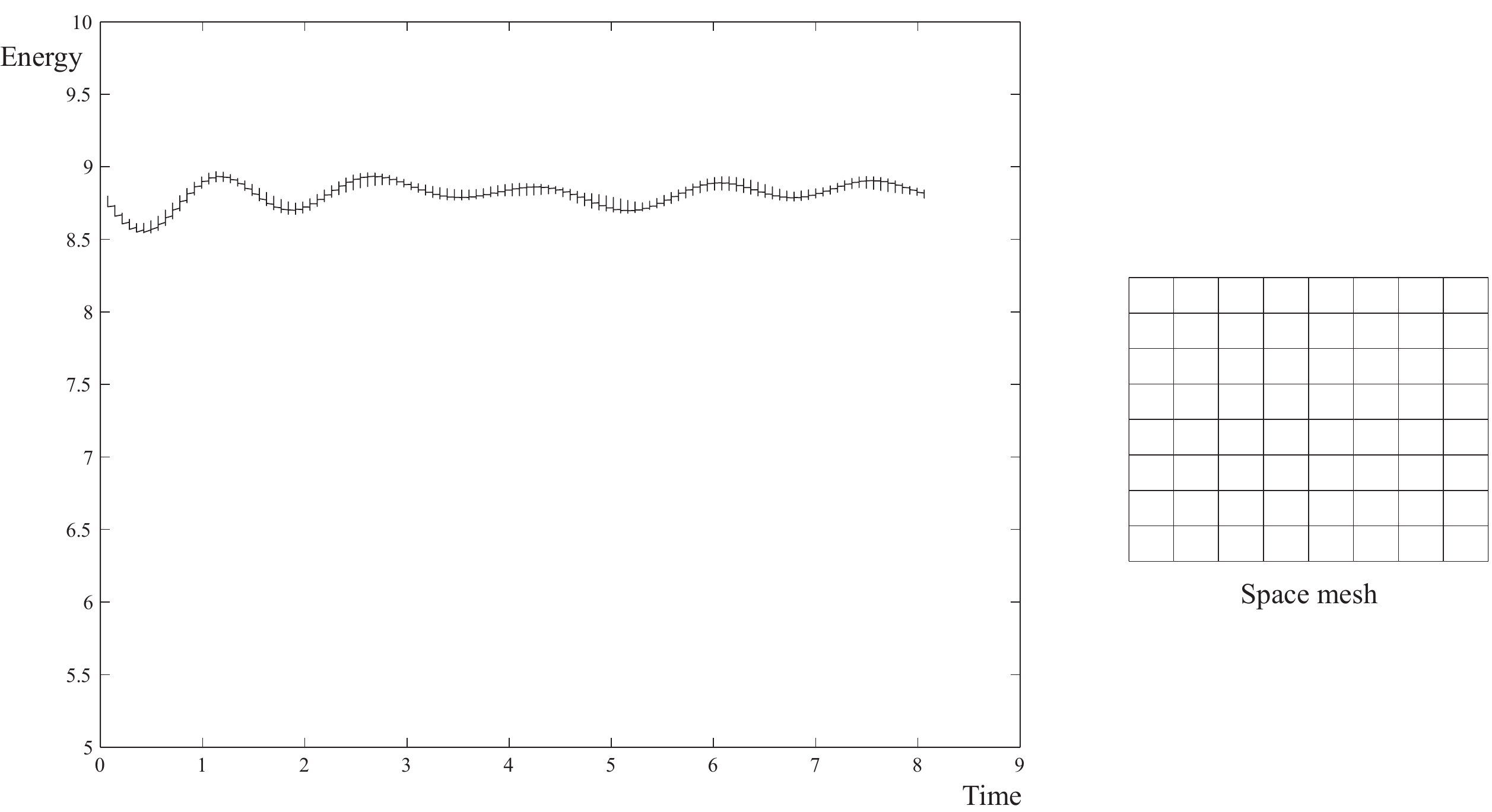}}
\caption[avi-energy-uniform-steps]{Energy vs. time for the AVI with
  uniform space and time discretization.  This is the special case
  where the AVI reproduces the Yee scheme --- which is well known to
  have good energy conservation properties, as seen here.  (The
  vertical ``tick marks'' on the plot show where the elements become
  synchronized, since they take uniform time
  steps.)} \label{fig:avi-energy-uniform-steps}
\end{figure}

\begin{figure}
\centerline{\includegraphics[width=0.7\linewidth]{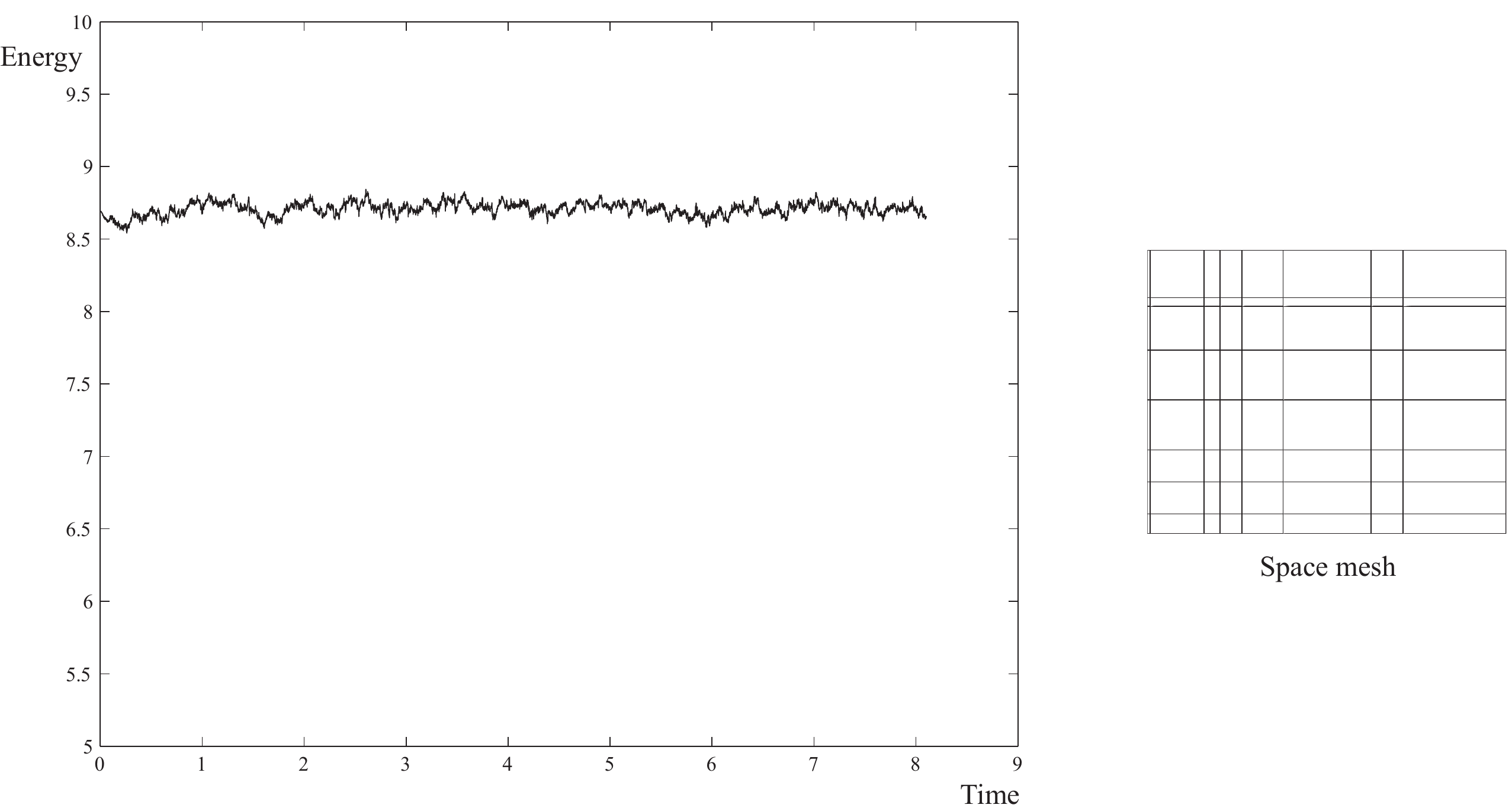}}
\caption[avi-energy-random-steps]{Energy vs. time for the AVI with
  random spatial discretization and fully asynchronous time steps.
  Despite the lack of regularity in the mesh and time steps, the AVI
  maintains the good energy behavior displayed by the Yee
  scheme.} \label{fig:avi-energy-random-steps}
\end{figure}

This simulation was done for two different spatial discretizations.
The first is a uniform discretization so that each element has
identical time step size, which coincides exactly with the Yee scheme.
The second discretization randomly partitioned the $x$- and $y$-axes,
so that each element has completely unique spatial dimensions and time
step size, and so the update rule is truly asynchronous.  The energy
plot for the uniform Yee discretization is shown in
\autoref{fig:avi-energy-uniform-steps}, while the energy for the
random discretization is shown in
\autoref{fig:avi-energy-random-steps}.  Even for a completely random,
irregular mesh, our asynchronous integrator displays near-energy
preservation qualities. Such numerical behavior stems from the
variational nature of our integrator, which will be detailed in
\autoref{section:theory}.

\begin{figure}
  \begin{center}
    \subfloat{\includegraphics[width=.6\linewidth]{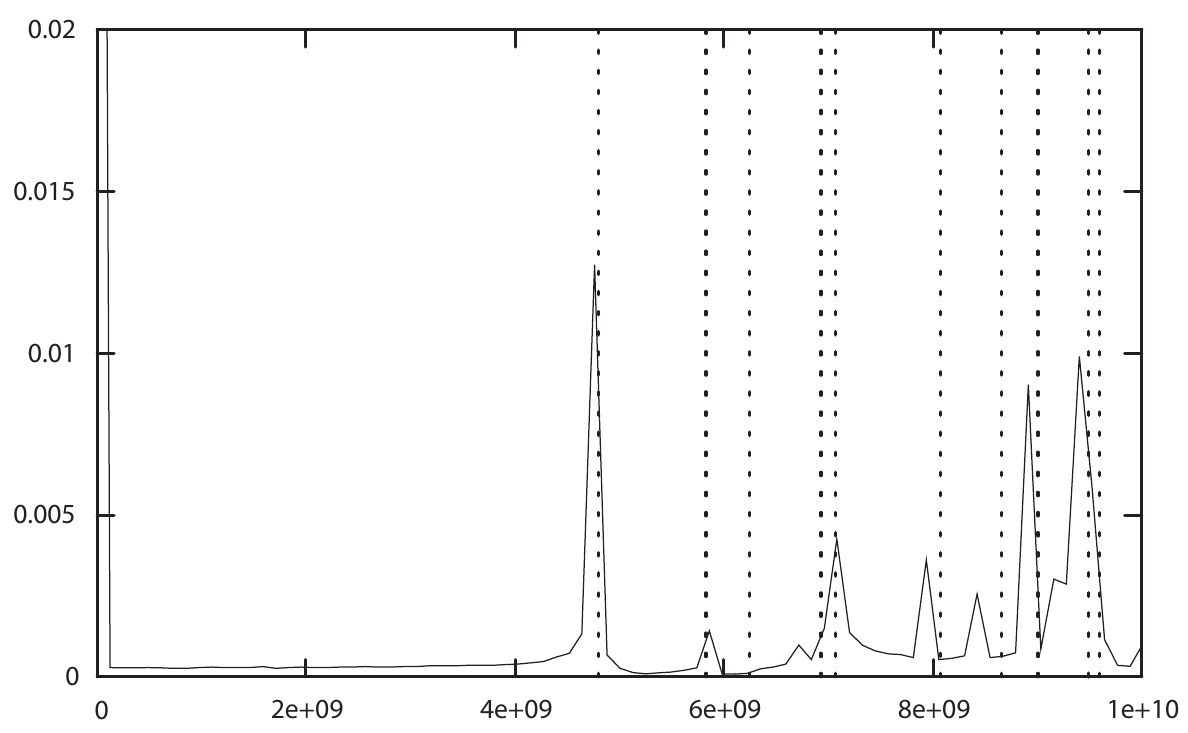}}
    \subfloat{\includegraphics[width=.35\linewidth]{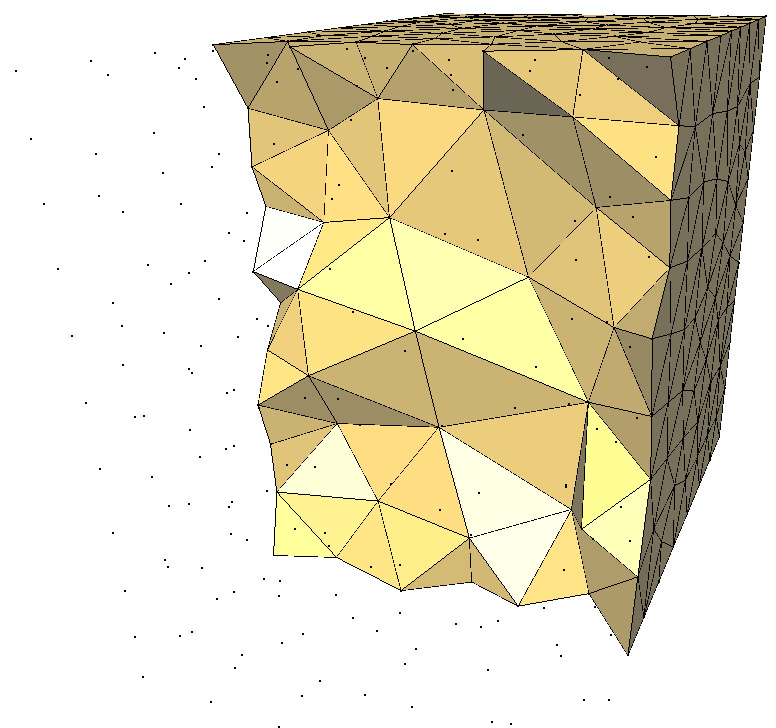}}
  \end{center}
  \caption{To produce the power spectrum shown at left, the electric
    field $E$ was initialized with random data (to excite all
    frequencies) and integrated forward in time, measuring the field
    strength at a particular sample point for every time step, and
    then performing a discrete Fourier transform.  The locations of
    the amplitude ``spikes'' are in close agreement with the analytic
    resonant frequencies, shown by the dashed vertical lines.  The
    spatial mesh, shown at right, was refined closer to the boundary,
    and coarser in the interior, allowing the AVI to produce this
    result with fewer total steps than uniform-time-stepping would
    require.}
  \label{fig:avi-resonance}
\end{figure}

In addition, we tested the performance of the AVI method with regard
to computing the resonant frequencies of a $3$-D rectangular cavity,
but using an {\em unstructured} tetrahedral spatial mesh.  While the
resonant frequencies are relatively simple to compute analytically,
nodal finite element methods are well known to produce spurious modes
for this type of simulation.  By contrast, as shown
in~\autoref{fig:avi-resonance}, the AVI simulation produces a
resonance spectrum in close agreement with theory.  Furthermore, by
refining the mesh close to the spatial boundary, while using a coarser
discretization in the interior, we were able to achieve these results
with less computational effort than a uniformly fine mesh would
require, since the time steps were selected to be proportional to the
respective element sizes.

\subsection{Fully Unstructured Spacetime Mesh}
\label{FullyUnstructured}

Finally, we look at the most general possible case: an arbitrary
discretization of spacetime, such as a simplicial $4$-complex.  Such a
mesh is completely relativistically covariant, so that $F$ cannot be
objectively separated into the components $E$ and $B$ without a
coordinate frame.  In most engineering applications, relativistic
effects are insignificant, so a 3+1 mesh (as in the previous
subsections) is almost always adequate, and avoids the additional
complications of spacetime mesh construction.  Still, we expect that
there are scientific applications where a covariant discretization of
electromagnetism may be very useful.  For example, many
implementations of numerical general relativity (using Regge calculus
for instance) are formulated on simplicial $4$-complexes; one might
wish to simulate the interaction of gravity with the electromagnetic
field, or charged matter, on such a mesh.

\subsubsection{Spacetime Mesh Construction}
First, a quick caution on mesh construction: since the Lorentz metric
is not positive definite, it is possible to create edges that have
length $0$, despite connecting two distinct points in $ \mathbb{R}
^{3,1} $ (so-called ``null'' or ``lightlike'' edges).  Meshes
containing such edges are degenerate---akin to a Euclidean mesh
containing a triangle with two identical points.  In particular, the
DEC Hodge star is undefined for $0$-volume elements (due to division
by zero).  Even without $0$-volume elements, it is still possible for
a spacetime mesh to violate causality, so extra care should be taken.
Methods to construct causality-respecting spacetime meshes over a
given spatial domain can be found in, e.g.,~\citet{ErGuSuUn2005}
and~\citet{Thite2005}.

When the mesh contains no inherent choice of a time direction, there
is no canonical way to split $F$ into $E$ and $B$.  Therefore, one
must set up the problem by assigning values of $F$ directly to
$2$-cells (or equivalently, assigning values of $A$ to $1$-cells).
For initial boundary value problems, one might choose to have the
initial and final time steps be prismal, so that $E$ and $B$ can be
used for initial and final values, while the internal discretization
is general.

\subsubsection{Equations of Motion}
The equations $ \mathrm{d} F = 0 $ and $ \mathrm{d} G = \mathcal{J} $
can be implemented directly in DEC.  Since this mesh is generally
unstructured, there is no simple algorithm as the ones we presented
above. Instead, the equations on $F$ results in a sparse linear system
which, given proper boundary conditions, can be solved globally with
direct or iterative solvers. However, it is clear that the previous
three examples that the methods of Yee, Bossavit--Kettunen, and our
AVI integrator are special cases where the global solution is
particularly simple to compute via synchronous or asynchronous time
updates.

\subsubsection{Mesh Construction and Energy Behavior}
It is known that, while variational integrators in mechanics do not
preserve energy exactly, they have excellent energy behavior, in that
it tends to oscillate close to the exact value.  This is only true,
however, when the integrator takes time steps of uniform size;
adaptive and other non-uniform stepping approaches can give poor
results unless additional measures are taken to enforce good energy
behavior.  (See \citealp[Chapter VIII]{HaLuWa2006} for a good
discussion of this problem for mechanics applications.)

Therefore, there is no reason to expect that \emph{arbitrary} meshes
of spacetime will yield energy results as good as the Yee,
Bossavit--Kettunen, and AVI schemes.  However, if one is taking a
truly covariant approach to spacetime, ``energy'' is not even
defined without specifying a time coordinate. Likewise, one would
not necessarily expect good energy behavior from the other methods
with respect to an arbitrary transformation of spatial coordinates.
Which sort of mesh to choose is thus highly application-dependent.

\section{Theoretical results}
\label{section:theory}

In this section, we complete our exposition with a number of
theoretical results about the discrete and continuous Maxwell's
equations.  In particular, we show that the DEC formulation of
electromagnetism derives from a discrete Lagrangian variational
principle, and that this formulation is consequently
multisymplectic. Furthermore, we explore the gauge symmetry of
Maxwell's equations, and detail how a particular choice of gauge
eliminates the equation for $ \nabla \cdot \mathbf{D} - \rho $ from
the Euler--Lagrange equations, while preserving it automatically as a
momentum map.

\begin{theorem}
The discrete Maxwell's equations are variational.
\end{theorem}

\begin{proof}
  The idea of this proof is to emulate the derivation of the
  continuous Maxwell's equations from \autoref{MaxwellEquations}.
  Interpreting this in the sense of DEC, we will obtain the discrete
  Maxwell's equations.

  Given a discrete $1$-form $A$ and dual source $3$-form
  $\mathcal{J}$, define the discrete Lagrangian $4$-form
\begin{equation}
  \mathcal{L} _d = -\frac{1}{2} \mathrm{d} A \wedge *\mathrm{d} A + A
  \wedge \mathcal{J} , \nonumber 
\end{equation}
with the corresponding discrete action principle
\begin{equation}
  S _d [A] = \left\langle \mathcal{L} _d , K \right\rangle . \nonumber 
\end{equation}
Then, taking a discrete $1$-form variation $\alpha$ vanishing on the
boundary, the corresponding variation of the action is
\begin{equation}
  \mathbf{d} S _d [A] \cdot \alpha = \left\langle - \mathrm{d} \alpha
    \wedge *\mathrm{d} A + \alpha \wedge \mathcal{J} , K \right\rangle =
  \left\langle \alpha \wedge \left( - \mathrm{d} {*\mathrm{d} A} +
      \mathcal{J} \right) , K \right\rangle  . \nonumber 
\end{equation}
(Here we use the bold $\mathbf{d}$ to indicate that we are
differentiating over the \emph{smooth space of discrete forms} $A$, as
opposed to differentiating over discrete spacetime, for which we use
$\mathrm{d}$.)  Setting this equal to $0$ for all variations $\alpha$,
the resulting discrete Euler--Lagrange equations are therefore $
\mathrm{d} {*\mathrm{d} A} = \mathcal{J} $.  Defining the discrete
$2$-forms $ F = \mathrm{d} A $ and $ G = *F $, this implies $
\mathrm{d} F = 0 $ and $ \mathrm{d}G = \mathcal{J} $, the discrete
Maxwell's equations.
\end{proof}

\subsection{Multisymplecticity}

The concept of {\em multisymplecticity} for Lagrangian field theories
was developed in~\citet{MaPaSh1998}, where it was shown to arise from
the boundary terms for general variations of the action, i.e., those
not restricted to vanish at the boundary.  As originally presented,
the Cartan form $\theta _{ \mathcal{L} } $ is an $ (n+1) $-form, where
the $n$-dimensional boundary integral is then obtained by contracting
$\theta _{ \mathcal{L} } $ with a variation.  The multisymplectic $
(n+2) $-form $\omega _{ \mathcal{L} } $ is then given by $ \omega _{
  \mathcal{L} } = - \mathbf{d} \theta _{ \mathcal{L} } $.  Contracting
$\omega _{ \mathcal{L} } $ with two arbitrary variations gives an
$n$-form that vanishes when integrated over the boundary, a result
called the {\em multisymplectic form formula}, which results from the
identity $ \mathbf{d} ^2 = 0 $.  In the special case of mechanics,
where $ n = 0 $, the boundary consists of the initial and final time
points; hence, this implies the usual result that the symplectic
$2$-form $\omega _L $ is preserved by the time flow.

\begin{figure}
  \centerline{\includegraphics[width=0.5\linewidth]{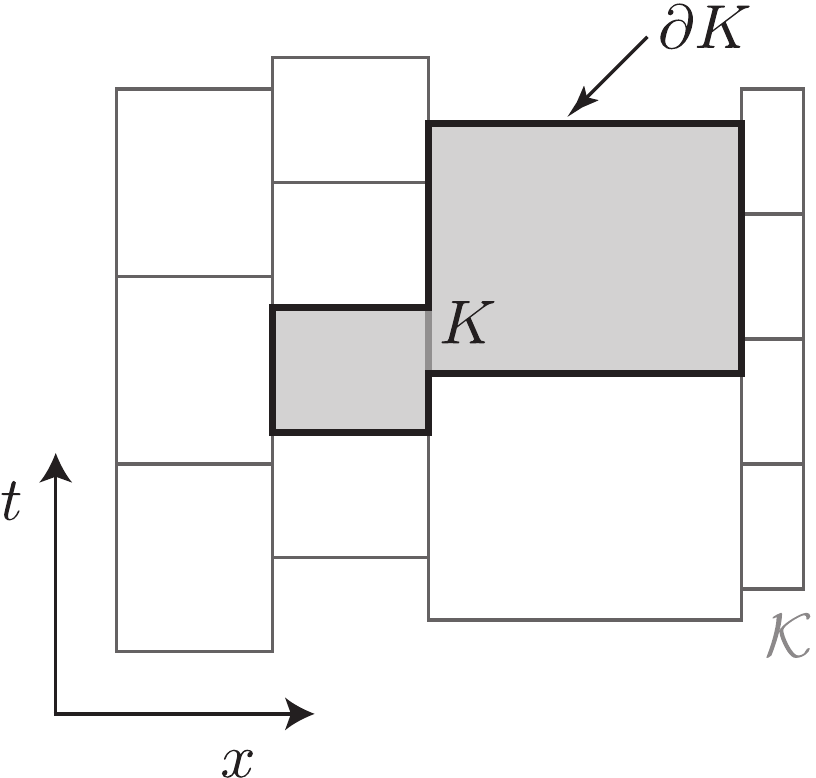}}
  \caption[multisymplectic]{To illustrate the discrete multisymplectic
    form formula \eqref{eq:multisymplectic}, we have here a $2$-{D}
    asynchronous-time mesh $\mathcal{K}$, where the shaded region is
    an arbitrary subcomplex $ K \subset \mathcal{K} $.  Given any two
    variations $ \alpha, \beta $ of the field, and the multisymplectic
    form $\omega _{ \mathcal{L} _d } $, the formula states that $
    \omega _{ \mathcal{L} _d } \cdot \alpha \cdot \beta $ vanishes
    when integrated over the boundary $ \partial K $ (shown in
    bold).} \label{fig:multisymplectic}
\end{figure}

Alternatively, as communicated to us by~\citet{Patrick2004}, one can
view the Cartan form $\theta _{\mathcal{L}} $ as an $n$-form-valued
$1$-form, and the multisymplectic form $\omega _{ \mathcal{L} } $ as
an $n$-form-valued $2$-form.  Therefore, one simply evaluates these
forms on tangent variations to obtain a boundary integral, rather than
taking contractions.  These two formulations are equivalent on smooth
spaces.  However, we will adopt Patrick's latter definition, since it
is more easily adapted to problems on discrete meshes: $\theta _{
  \mathcal{L} } $ and $\omega _{ \mathcal{L} } $ remain smooth $1$-
and $2$-forms, respectively, but their $n$-form values are now taken
to be discrete.  See \autoref{fig:multisymplectic} for an illustration
of the discrete multisymplectic form formula.

\begin{theorem}
The discrete Maxwell's equations are multisymplectic.
\end{theorem}

\begin{proof}
  Let $ K \subset \mathcal{K} $ be an arbitrary subcomplex, and
  consider the discrete action functional $S _d $ restricted to $K$.
  Suppose now that we take a discrete variation $\alpha$, {\em
    without} requiring it to vanish on the boundary $\partial K$.
  Then variations of the action contain an additional boundary term
\begin{equation}
  \mathbf{d} S _d [A] \cdot \alpha = \left\langle \alpha \wedge \left(
      - \mathrm{d} {*\mathrm{d} A } + \mathcal{J} \right) , K \right\rangle 
  + \left\langle \alpha \wedge *d A, \partial K \right\rangle . \nonumber 
\end{equation}
Restricting to the space of potentials $A$ that satisfy the discrete
Euler--Lagrange equations, the first term vanishes, leaving only
\begin{equation}
  \mathbf{d} S _d (A) \cdot \alpha = \left\langle \alpha \wedge
    *\mathrm{d} A,\partial K \right\rangle \label{cartanform}
\end{equation}

Then we can define the \emph{Cartan form} $\theta _{ \mathcal{L} _d }
$ by
\begin{equation}
  \theta _{ \mathcal{L} _d } \cdot \alpha  = \alpha  \wedge
  *\mathrm{d} A. \nonumber 
\end{equation}
Since $\theta _{ \mathcal{L} _d } $ takes a tangent vector $\alpha$
and produces a discrete $3$-form on the boundary of the subcomplex, it
is a {\em smooth 1-form taking discrete $3$-form values}.  Now, since
the {\em space of discrete forms} is itself actually continuous, we
can take the exterior derivative in the smooth sense on both sides of
\autoref{cartanform}.  Evaluating along another first variation
$\beta$ (again restricted to the space of Euler--Lagrange solutions),
we then get
\begin{equation}
  \mathbf{d} ^2 S _d [A] \cdot \alpha \cdot \beta = \left\langle
    \mathbf{d} \theta \cdot \alpha \cdot \beta , \partial K
  \right\rangle. \nonumber 
\end{equation}
Finally, defining the multisymplectic form $ \omega _{ \mathcal{L}
  _{d} } = - \mathbf{d} \theta _{ \mathcal{L} _d } $, and using the
fact that $ \mathbf{d} ^2 S_{d} = 0 $, we get the relation
\begin{equation}
  \left\langle \omega _{ \mathcal{L} _d } \cdot \alpha \cdot \beta , \partial K
  \right\rangle = 0
  \label{eq:multisymplectic}
\end{equation}
for all variations $ \alpha , \beta$;
\autoref{eq:multisymplectic} is a discrete version of the
multisymplectic form formula.  Since this holds for any subcomplex
$K$, it follows that these schemes are multisymplectic.
\end{proof}

\subsection{Gauge Symmetry Reduction and Covariant Momentum Maps}
\label{section:gauge}

We now explore the symmetry of Maxwell's equations under gauge
transformations.  This symmetry allows us to reduce the equations by
eliminating the time component of $A$ (for some chosen time
coordinate), effectively fixing the electric scalar potential to zero.
Because this is an incomplete gauge, there is a remaining gauge
symmetry, and hence a conserved momentum map.  This conserved quantity
turns out to be the charge density $ \rho = \nabla \cdot \mathbf{D} $,
which justifies its elimination from the Euler--Lagrange equations.
These calculations are done with differential forms and exterior
calculus, hence they apply equally to the smooth and discrete cases of
electromagnetism.

\subsection{Choosing a Gauge}
Because Maxwell's equations only depend on $ \mathrm{d} A $, they are
invariant under gauge transformations of the form $ A \mapsto A +
\mathrm{d} f $, where $f$ is any scalar function on spacetime.  If we
fix a time coordinate, we can now choose the Weyl gauge, so that the
time component $ A _t = 0 $.  Therefore, we can assume that
\begin{equation}
  A = A _x \;\mathrm{d}x  + A _y \;\mathrm{d}y  + A _z \;\mathrm{d}z  . \nonumber 
\end{equation}
In fact, $ A _x, A _y, A _z $ are precisely the components of the
familiar vector potential $\mathbf{A}$, i.e., $ A = \mathbf{A} ^\flat
$.

\subsection{Reducing the Equations}
Having fixed the gauge and chosen a time coordinate, we can now define
two new ``partial exterior derivative'' operators, $ \mathrm{d} _t $
(time) and $ \mathrm{d} _s $ (space), where $ \mathrm{d} = \mathrm{d}
_t + \mathrm{d} _s $.  Since $A$ contains no $ \mathrm{d} t $ terms, $
\mathrm{d} _s A $ is a $2$-form containing only the space terms of $
\mathrm{d} A $, while $ \mathrm{d} _t A $ contains the terms involving
both space and time.  That is,
\begin{eqnarray}
  \mathrm{d}  _t A = E \wedge \mathrm{d}  t, && \mathrm{d} _s A = B. \nonumber 
\end{eqnarray}
Restricted to this subspace of potentials, the Lagrangian density then
becomes
\begin{align*}
  \mathcal{L} &= -\frac{1}{2} \left( \mathrm{d} _t A + \mathrm{d} _s A
  \right) \wedge
  *\left( \mathrm{d}  _t A + \mathrm{d}  _s A \right)  + A \wedge \mathcal{J} \\
  &= -\frac{1}{2} \left( \mathrm{d} _t A \wedge * \mathrm{d} _t A +
    \mathrm{d} _s A \wedge * \mathrm{d} _s A \right) + A \wedge J
  \wedge \mathrm{d} t
\end{align*}
Next, varying the action along a restricted variation $\alpha$ that
vanishes on $ \partial X $,
\begin{align}
  \mathbf{d} S[A] \cdot \alpha &= \int _X \left( \mathrm{d} _t \alpha
    \wedge {D} - \mathrm{d} _s \alpha \wedge H \wedge
    \mathrm{d} t + \alpha \wedge J \wedge \mathrm{d} t
  \right) \label{eqn:weylaction}\\
  &= \int _X \alpha \wedge \left( \mathrm{d} _t D - \mathrm{d} _s H
    \wedge \mathrm{d} t + J \wedge \mathrm{d} t \right) .\nonumber
\end{align}
Setting this equal to zero by Hamilton's principle, one immediately
gets Amp\`ere's law as the sole Euler--Lagrange equation.  The
divergence constraint $ \mathrm{d} _s {D} = \rho $, corresponding to
Gauss' law, has been eliminated via the restriction to the Weyl gauge.

\subsubsection{Noether's Theorem Implies Automatic Preservation of
  Gauss' Law}
Let us restrict $A$ to be an Euler--Lagrange solution in the Weyl
gauge, but remove the previous requirement that variations $\alpha$ be
fixed at the initial time $ t _0 $ and final time $ t _f $.  Then,
varying the action along this new $\alpha$, the Euler--Lagrange term
disappears, but we now pick up an additional boundary term due to
integration by parts
\begin{equation*}
  \mathbf{d} S [A] \cdot \alpha = \left. \int _{\Sigma} \alpha \wedge
    {D} \,\right\rvert _{t _0 } ^{ t _f } ,
\end{equation*}
where $\Sigma$ denotes a Cauchy surface of $X$, corresponding to the
spatial domain.  If we vary along a gauge transformation $ \alpha =
\mathrm{d} _s f $, then this becomes
\begin{equation*}
  \mathbf{d} S [A] \cdot \mathrm{d}  _s f = \left. \int _\Sigma
    \mathrm{d}  _s f \wedge {D} \,\right\rvert _{t _0 } ^{ t _f }
  = - \left. \int _\Sigma f \wedge \mathrm{d} _s {D} \,\right\rvert _{t _0 }
  ^{ t _f }
\end{equation*}
Alternatively, plugging $ \alpha = \mathrm{d} _s f $ into
\autoref{eqn:weylaction}, we get
\begin{equation*}
  \mathbf{d} S [A] \cdot \mathrm{d}  _s f = \int _X \mathrm{d} _s
  f \wedge J \wedge \mathrm{d}  t
  = - \int _X f \wedge \mathrm{d}  _s J \wedge \mathrm{d}  t
  = - \int _X f \wedge \mathrm{d} _t \rho
  = - \left. \int _\Sigma f \wedge \rho \,\right\rvert _{ t _0 } ^{ t _f }
  .
\end{equation*}
Since these two expressions are equal, and $f$ is an arbitrary
function, it follows that
\begin{equation*}
  \left. \left( \mathrm{d}  _s D - \rho \right) \right\rvert _{t_0}^{t_f} = 0.
\end{equation*}
This indicates that $ \mathrm{d} _s D - \rho $ is a conserved
quantity, a momentum map, so if Gauss' law holds at the initial time,
then it holds for all subsequent times as well.

\subsection{Boundary Conditions and Variational Structure}
It should be noted that the variational structure and symmetry of
Maxwell's equations may be affected by the boundary conditions that
one chooses to impose.  There are many boundary conditions that one
can specify independent of the initial values, such as the PEC
condition used in the numerical example in \autoref{AVI}.  However,
one can imagine more complicated boundary conditions where which the
boundary interacts nontrivially with the interior of the domain ---
such as dissipative or forced boundary conditions, where
energy/momentum is removed from or added to the system.  In these
cases, one will obviously {\em not} conclude that the charge density $
\nabla \cdot \mathbf{D} $ is conserved, but more generally that the
{\em change} in charge is related to the flux through the spatial
boundary.  This is because, in the momentum map derivation above, the
values of $f$ on the initial time slice causally affects its values on
the spatial boundary at intermediate times, not just on the final time
slice.  Thus, the spatial part of $ \partial X $ cannot be neglected
for arbitrary boundary conditions.

\section{Conclusion}

The continued success of the Yee scheme for many applications of
computational electromagnetism, for over four decades, illustrates the
value of structure-preserving numerical integrators for Maxwell's
equations.  Recent advances by, among others, Bossavit and Kettunen,
and Gross and Kotiuga, have demonstrated the important role of
compatible spatial discretization using differential forms, allowing
for Yee-like schemes that apply on generalized spatial meshes.  In
this paper, we have extended this approach by considering discrete
forms on {\em spacetime}, encapsulating both space and time
discretization, and have derived a general family of geometric
numerical integrators for Maxwell's equations.  Furthermore, since we
have derived these integrators from a discrete variational principle,
the resulting methods are provably multisymplectic and
momentum-map-preserving, and they experimentally show correct global
energy behavior.  Besides proving the variational nature of well-known
techniques such as the Yee and Bossavit--Kettunen schemes, we have also
introduced a new asynchronous integrator, so that time step sizes can
be taken non-uniformly over the spatial domain for increased
efficiency, while still maintaining the desirable variational and
energy behavior of the other methods.

\subsection*{Future Work}
One promising avenue for future work involves increasing the order of
accuracy of these methods by deriving higher-order discrete Hodge star
operators.  While this would involve redefining the Hodge star matrix
to be non-diagonal, the discrete Maxwell's equations would remain
formally the same, and hence there would be no change in the
variational or multisymplectic properties proven here.  We are
currently exploring the development of a spectrally accurate spatial
Hodge star, which might make these geometric schemes competitive for
applications where non-variational spectral codes are currently
favored.

Additionally, the recent work of~\citet{KaLe2007} has shown that AVIs
can be implemented as parallel algorithms for solid mechanics
simulations.  This uses the fact that, due to the asynchronous update
procedure, an element does not need information from every one of its
neighbors at every time step, which lessens the need for communication
among parallel nodes.  The resulting parallel AVIs, or PAVIs, can
therefore take advantage of parallel computing architecture for
improved efficiency.  It is reasonable to expect that the same might
be done in the case of our electromagnetic AVI.

While we have experimentally observed the fact that variational
integrators exhibit near-energy conservation, little is known about
this behavior from a theoretical standpoint.  In the case of ODEs in
mechanics, backwards error analysis has shown that these methods
exactly integrate a nearby smooth Hamiltonian system, although not
much known about how this relates to the discrete variational
principle on the Lagrangian side.  Some initial work has been done
in~\citet{OlWeWu2004} to understand, also by a backward error analysis
approach, why discrete multisymplectic methods also display good
energy behavior.

Finally, variational methods using discrete spacetime forms may be
developed for field theories other than electromagnetism.  Promising
candidates include numerical general relativity and fluid dynamics,
although the latter is complicated by the difficulty in finding a
proper discretization of the infinite-dimensional diffeomorphism
group.  If discrete Lagrangian densities are developed for these
theories, it should be straightforward to combine them with the
electromagnetic Lagrangian, resulting in numerical methods to
simulate, e.g., gravity coupled with an electromagnetic field, or the
dynamics of a charged or magnetic fluid.

\section*{Acknowledgments}
We would like to thank several people for their inspiration and
suggestions. First of all, Alain Bossavit for suggesting many years
ago that we take the present DEC approach to computational
electromagnetism, and for his excellent lectures at Caltech on the
subject. Second, Michael Ortiz and Eva Kanso for their ongoing
interactions on related topics and suggestions.  We also thank Doug
Arnold, Uri Ascher, Robert Kotiuga, Melvin Leok, Adrian Lew, and Matt
West for their feedback and encouragement.  In addition, the $3$-{D}
AVI simulations shown in \autoref{fig:avi-resonance} were programmed
and implemented by Patrick Xia, as part of a Summer Undergraduate
Research Fellowship at Caltech.

\end{document}